\documentclass[11pt]{article}
\def\revisiondate{20 December 2013}

\usepackage{graphicx}
\usepackage{amssymb}
\usepackage{amsmath}
\usepackage{amsthm}
\usepackage{afterpage}
\usepackage{hyperref}

\hypersetup{%
  bookmarksopen=true,
  bookmarksopenlevel=2,
  bookmarksnumbered=true,
  pdftitle={Non-Embeddable Extensions of Embedded Minors},
  pdfauthor={Rajneesh Hegde},
  pdfkeywords={Graph Embeddings, Non-planar Graphs, Graph Minors},
  pageanchor=true,
  colorlinks=true,
  linkcolor={blue},
  citecolor={blue},
  urlcolor={red},
  pagecolor={cyan},
  anchorcolor={black}
}

\newcommand{\wfc}{weakly $4$-connected}
\newcommand{\wfcity}{weak $4$-connectivity}

\newcommand{\gp}{\ensuremath{G^+}}
\newcommand{\gpp}{\ensuremath{G^{\scriptscriptstyle ++}}}
\newcommand{\ghat}{\ensuremath{\widehat{G}}}

\newcommand{\D}{\ensuremath{\mathcal{D}}}

\newcommand{\Dat}{\ensuremath{\widetilde{D}_1}}
\newcommand{\Dbt}{\ensuremath{\widetilde{D}_2}}
\newcommand{\ext}[1]{{\rm#1}-enlargement}

\newcommand{\jump}{\ref{jump}}
\newcommand{\cross}{\ref{cross}}
\newcommand{\ncsplit}{\ref{ncsplit}}
\newcommand{\splitjump}{\ref{split+jump}}
\newcommand{\dsplitjump}{\ref{dsplit+jump}}
\newcommand{\splitcross}{\ref{split+cross}}
\newcommand{\dsplitcross}{\ref{dsplit+cross}}
\newcommand{\triad}{\ref{triad}}
\newcommand{\tedge}{\ref{t-edge}}

\newcommand{\prism}{\ref{prism}}

\newcommand{\fsplit}{\ensuremath{F_1'}}
\newcommand{\ffour}{\ensuremath{F_4}}
\newcommand{\dsplit}{\ensuremath{D_3'}}
\newcommand{\kffmfour}{\ensuremath{E_{22}}}
\newcommand{\etwenty}{\ensuremath{E_{20}}}
\newcommand{\cthree}{\ensuremath{C_3}}
\newcommand{\etwo}{\ensuremath{E_2}}
\newcommand{\kffminus}{\ensuremath{E_{18}}}

\newtheorem{lemma}{Lemma}[section]
\newtheorem{theorem}[lemma]{Theorem}

\newcommand{\figref}[1]{\hyperref[#1]{Figure~\ref*{#1}}}%
\newcommand{\tabref}[1]{\hyperref[#1]{Table~\ref*{#1}}}%
\newcommand{\secref}[1]{\hyperref[#1]{Section~\ref*{#1}}}%
\newcommand{\chapref}[1]{\hyperref[#1]{Chapter~\ref*{#1}}}%
\newcommand{\thmref}[1]{\hyperref[#1]{Theorem~\ref*{#1}}}%
\newcommand{\lemref}[1]{\hyperref[#1]{Lemma~\ref*{#1}}}%
\newcommand{\obsref}[1]{\hyperref[#1]{Observation~\ref*{#1}}}%
\newcommand{\corref}[1]{\hyperref[#1]{Corollary~\ref*{#1}}}%
\newcommand{\conjref}[1]{\hyperref[#1]{Conjecture~\ref*{#1}}}%

\newenvironment{myenumerate}[2][\arabic]{\begin{enumerate}%
\renewcommand{\theenumi}{(#2#1{enumi})}%
\renewcommand{\labelenumi}{\theenumi}}%
{\end{enumerate}}

\title{Non-Embeddable Extensions of Embedded Minors%
  \thanks{Partially supported by NSF grants DMS-0200595 and 0354742.}
}
\author{Rajneesh Hegde and Robin Thomas}

\begin{document}
\phantom{a}\vskip .25in
\centerline{\Large{\bf NON-EMBEDDABLE EXTENSIONS OF EMBEDDED MINORS}%
  \footnote{Partially supported by NSF grants DMS-0200595 and 0354742.
  \revisiondate.}
}
\bigskip
\centerline{{\bf Rajneesh Hegde} and {\bf Robin Thomas}}
\centerline{School of Mathematics}
\centerline{Georgia Institute of Technology}
\centerline{Atlanta, Georgia  30332, USA}
\bigskip

\begin{abstract}
  A graph $G$ is \wfc\ if it is $3$-connected, has at least five
  vertices, and for every pair $(A,B)$ such that $A\cup B=V(G)$,
  $|A\cap B|=3$ and no edge has one end in $A-B$ and the other in
  $B-A$, one of the induced subgraphs $G[A], G[B]$ has at most four
  edges. We describe a set of constructions that starting from a
  \wfc\ planar graph $G$ produce a finite list of non-planar \wfc\
  graphs, each having a minor isomorphic to $G$, such that every
  non-planar \wfc\ graph $H$ that has a minor isomorphic to $G$ has a
  minor isomorphic to one of the graphs in the list. Our main result
  is more general and applies in particular to polyhedral embeddings
  in any surface.
\end{abstract}

\baselineskip 18pt

\section{Introduction}
\label{sec:disks-intro}

We begin with some basic notation and ingredients needed to state the
main result of this paper. Graphs are finite and simple (i.e., they
have no loops or multiple edges). \emph{Paths} and \emph{cycles} have
no ``repeated" vertices. A graph is a \emph{minor} of another if the
first can be obtained from a subgraph of the second by contracting
edges. For a graph $G$ and an edge $e$ in $G$, $G \backslash e$ an
$G/e$ are the graphs obtained from $G$ by respectively deleting and
contracting the edge $e$. A graph is a \emph{subdivision} of another
if the first can be obtained from the second by replacing each edge by
a non-zero length path with the same ends, where the paths are
disjoint, except possibly for shared ends. The replacement paths are
called \emph{segments}, and their ends are called
\emph{branch-vertices}. A graph is a \emph{topological minor} of
another if a subdivision of the first is a subgraph of the second.

Let a non-planar graph $H$ have a subgraph isomorphic to a subdivision
of a planar graph $G$. For various problems in Graph Structure Theory
it is useful to know the minimal subgraphs of $H$ that are isomorphic
to a subdivision of $G$ and are non-planar. In other words, one wants
to know what more does $H$ contain on account of its non-planarity.
In~\cite{RobSeyThoExt} it is shown that under some mild connectivity
assumptions these ``minimal non-planar enlargements'' of $G$ are quite
nice. In the applications of the result, $G$ is explicitly known,
whereas $H$ is not, and the enlargement operations would furnish an
explicit list of graphs such that (i)~$H$ has a subgraph isomorphic to
a subdivision of one of the graphs on the list, and (ii)~each graph on
the list is a witness both to the fact that $G$ is a topological minor
of $H$, and that $H$ is, in addition, non-planar. (The minimality of
the graphs in the list is required to avoid redundancy.) Before we
state that result, we need a few definitions.

For $Z \subseteq V(G)$, $G[Z]$ denotes the subgraph \emph{induced} by
$Z$, that is, the subgraph consisting of $Z$ and all edges with both
ends in $Z$. A subgraph of $G$ is said to be \emph{induced} if it is
induced by its vertex set.

A \emph{separation} of a graph $G$ is a pair $(A,B)$ of subsets of
$V(G)$ such that $A \cup B = V(G)$, and there is no edge between $A-B$
and $B-A$. The \emph{order} of $(A,B)$ is $|A \cap B|$. The separation
is called \emph{non-trivial} if both $A$ and $B$ are proper subsets of
$V(G)$.
A graph $G$ is \emph{weakly $4$-connected} if $G$ is
$3$-connected, has at least five vertices, and for every separation
$(A,B)$ of $G$ of order at most three, one of the graphs $G[A]$, $G[B]$
has at most four edges. 

A cycle $C$ in a graph $G$ is called \emph{peripheral} if it is
induced and $G\backslash V(C)$ is connected. It is
well-known~\cite{TutHowto,WhiCongr} that the peripheral cycles in a
$3$-connected planar graph are precisely the cycles that bound faces
in some (or, equivalently, every) planar embedding of $G$.

Let $S$ be a subgraph of a graph $H$. An \emph{$S$-path} in $H$ is a path
with both ends in $S$, and otherwise disjoint from $S$.
Let $C$ be a cycle in  $S$, and let $P_1$ and $P_2$ be two
disjoint $S$-paths in $H$ with ends $u_1,v_1$ and $u_2,v_2$, respectively,
such that $u_1,u_2,v_1,v_2$ belong to $V(C)$ and occur on $C$ in the
order listed. In those circumstances we say that the pair $P_1,P_2$ is
an \emph{$S$-cross in $H$}. We also say that it is an \emph{$S$-cross on
  $C$.} We say that $u_1,v_1,u_2,v_2$ are the \emph{feet} of the
cross. We say that the cross $P_1,P_2$ is \emph{free} if
\begin{myenumerate}{F}
\item for $i=1,2$ no segment of $S$ includes both ends of $P_i$,
  and
\item no two segments of $S$ that share a vertex include all the
  feet of the cross.
\end{myenumerate}
The following was proved in~\cite{RobSeyThoExt}.

\begin{theorem}
\label{thm:top-nonpl-ext}
Let $G$ be a weakly $4$-connected planar graph, 
and let $H$ be a weakly $4$-connected non-planar graph such that
a subdivision of $G$ is isomorphic to a subgraph of $H$. Then
there exists a subgraph $S$ of $H$ isomorphic to a subdivision of
$G$ such that one of the following conditions holds:
\begin{enumerate}
\item there exists an $S$-path in $H$ such that its ends belong to no
  common peripheral cycle in $S$, or
\item there exists a free $S$-cross in $H$ on some peripheral cycle of
  $S$.
\end{enumerate}
\end{theorem}

This theorem has been used in \cite{DinOpoThoVer,ThoThoTutte}, and
its extension has been used in~\cite{KawNorThoWolbdtw}.
However, in more
complicated applications it is more efficient to work with minors,
rather than topological minors. We sketch one such application in
\secref{sec:application}.
For any fixed graph $G$, there exists a finite and
explicitly constructible set $\{G_1,G_2,\ldots,G_t\}$ of graphs such
that a graph $H$ has a minor isomorphic to $G$ if and only if it has a
topological minor isomorphic to one of the graphs $G_i$. 
Thus one can apply \thmref{thm:top-nonpl-ext} $t$ times to deduce
the desired conclusion about $G$, but it would be nicer to have a more
direct route to the result that involves less potential
duplication. Furthermore, if the outcome is allowed to be a minor of
$H$ rather than a topological minor, then the outcomes~(i) and (ii)
above can be strengthened to require that the ends of the paths
involved are branch-vertices of $S$, as we shall see.

It turns out that \thmref{thm:top-nonpl-ext} is not
exclusively about face boundaries of planar graphs, but that an
appropriate generalization holds under more general circumstances.
Thus rather than working with peripheral cycles in planar graphs we
will introduce an appropriate set of axioms for a set of cycles of a
general graph. We do so now in order to avoid having to restate our
definitions later when we present the more general form of our
results.

A \emph{segment} in a graph $G$ is a maximal path 
such that its internal
vertices all have degree in $G$ exactly two. 
If a graph $G$ has no vertices of degree two, then the segments of
a subdivision of $G$ defined earlier coincide with the notion just
defined. Since we will not consider subdivisions of graphs with 
vertices of degree two there is no danger of confusion.
A \emph{cycle double cover} in a graph $G$ is a set \D\ of distinct cycles
of $G$, called \emph{disks}, such that
\newcounter{axctr}
\begin{myenumerate}{D}
\item \label{ax:1} each edge of $G$ belongs to precisely two members of \D.
\setcounter{axctr}{\value{enumi}}
\end{myenumerate}
A cycle double cover \D\ is called
a \emph{disk system} in $G$ if 
\begin{myenumerate}{D}
\setcounter{enumi}{\value{axctr}}
\item \label{ax:2} for every vertex $v$ of $G$, the edges
  incident with $v$ can be arranged in a cyclic order such that for
  every pair of consecutive edges in this order, there is precisely one
  disk in \D\ containing that pair of edges, and
\item \label{ax:3}  the intersection of any two distinct disks in \D\ either
  has at most one vertex or is a segment.
\end{myenumerate}

A cycle double cover satisfying~\ref{ax:3} is called a \emph{weak disk
  system}.  It is easy to see that if a connected graph has a disk
system, then it is a subdivision of a $3$-connected graph. Also, note
that in a 3-connected graph, Axiom~\ref{ax:3} is equivalent to the
requirement that every two distinct disks intersect in a complete
subgraph on at most two vertices. The peripheral cycles of a
$3$-connected planar graph form a disk system. More generally, if $G$
is a subdivision of a $3$-connected graph embedded in a surface
$\Sigma$ in such a way that every homotopically nontrivial closed
curve intersects the graph at least three times (a ``polyhedral
embedding''), then the face boundaries of this embedding form a disk
system in $G$. Conversely, it can be shown that a disk system in a
graph is the set of face boundaries of a polyhedral embedding of the
graph in some surface.  Weak disk systems correspond to face
boundaries of embeddings into pseudosurfaces (surfaces with ``pinched"
points).

Let $G$ be a graph with a cycle double cover \D. Two vertices or
edges of $G$ are said to be \emph{confluent} if there is a disk
containing both of them. If \D\ is a cycle double cover in a graph
$G$ and $S$ is a subdivision of $G$, then \D\ induces a cycle double
cover $\D'$ in $S$ in the obvious way, and vice versa. We say that
$\D'$ is the \emph{cycle double cover induced in $S$ by \D}.


Let $v$ be a vertex of a graph $G$ with degree at least 4. Partition
the set of its neighbors into two disjoint sets $N_1$ and $N_2$, with
at least two vertices in each set. Let $G'$ be obtained from $G$ by
replacing the vertex $v$ with two adjacent vertices $v_1,v_2$, with
$v_i$ adjacent to the vertices in $N_i$ for $i=1,2$. The graph $G'$ is
said to be obtained from $G$ by \emph{splitting} the vertex $v$. It is
easy to see that if $G$ is 3-connected, then so is $G'$. The vertices
$v_1$ and $v_2$ are called the \emph{new vertices} of $G'$ and the
edge $v_1v_2$ of $G'$ is called the \emph{new edge} of $G'$.

Suppose a graph $G$ has a cycle double cover \D. The above
splitting operation on a vertex $v$ of $G$ is said to be a
\emph{conforming split} (with respect to \D) if
\begin{myenumerate}{S}
\item \label{split:1} among the disks that use the vertex $v$, there
  are exactly two, say $D_1$ and $D_2$, that use one vertex each from
  $N_1$ and $N_2$, and
\item \label{split:2} $D_1$ and $D_2$ intersect precisely in the vertex
  $v$.
\end{myenumerate}
The split is then said to be \emph{along} $D_1$ (and along $D_2$). A
split that is not conforming as above is said to be a
\emph{non-conforming split}.

Let $G,G'$ be as in the above paragraph.  If $G$ is a 3-connected
planar graph, then $G'$ is planar if and only if the split is
conforming with respect to the disk system of peripheral cycles of
$G$. More generally, to each cycle $C$ of $G$ there corresponds a
unique cycle $C'$ of $G'$, and so to \D\ there corresponds a uniquely
defined set of cycles $\D'$ of $G'$. If \D\ is a weak disk system,
then so is $\D'$, and if \D\ is a disk system, then so is $\D'$. We
call $\D'$ the (weak) disk system \emph{induced} in $G'$ by \D. This
is the purpose of conditions~\ref{split:1} and \ref{split:2}. If \D\
is a disk system, then an equivalent way to define a conforming split
of a vertex $v$ is to say that both $N_1$ and $N_2$ form contiguous
intervals in the cyclic order induced on the neighborhood of $v$ by
\D. Similarly, an equivalent condition for a split to be
non-conforming with respect to a disk system is the existence of
vertices $a,c \in N_1$ and $b,d \in N_2$ such that $a,b,c$ and $d$
appear in the cyclic order listed around $v$ (as given by \D\ in
\ref{ax:2}). The reason we use the definition above is that it
applies more generally to weak disk systems.



A graph $G'$ obtained from a graph $G$ by repeatedly splitting
vertices of degree at least four is said to be an \emph{expansion} of
$G$. In particular, each graph is an expansion of itself. Each split
leading to an expansion of $G$ has exactly one new edge; the set of
these edges are the \emph{new edges of the expansion of $G$}. The new
edges form a forest in $G'$. If $G$ has a cycle double cover \D, the
expansion is called a \emph{conforming} expansion if each of the
splits involved in it is conforming (with respect to \D). If at least
one of the splits involved is not conforming, then the expansion is
called \emph{non-conforming}. From the above discussion, it is clear
that a disk system in $G$ induces a unique disk system in a conforming
expansion.

We now describe seven enlargement operations.  Let $G$ be a graph with
a cycle double cover \D, and let \gp\ be the graph obtained from $G$
by applying one of the operations described below.
\newcounter{opctr}
\begin{enumerate}
\item (non-conforming jump) \label{jump} \gp\ is obtained from $G$ by
  adding an edge $uv$ where $u$ and $v$ are non-confluent vertices of $G$.
\item (cross) \label{cross} Let $a,b,c,d$ be vertices appearing on a
  disk of $G$ in that cyclic order. Add the edges $ac$ and $bd$
  to obtain \gp.
\item (non-conforming split) \label{ncsplit} \gp\ is obtained from $G$ by
  performing a non-conforming split of a vertex of $G$.
\item (split + non-conforming jump) \label{split+jump} Let $u,v$ be
  non-adjacent vertices on some disk $C \in \D$. Perform a conforming
  split of $v$ into $v_1,v_2$ such that $u$ and $v_2$ are
  non-confluent vertices.  (In particular, the split is not along $C$.)
  Now add the edge $uv_2$ to obtain \gp.
\item (double split + non-conforming jump) \label{dsplit+jump} Let $u,v$
  be adjacent vertices and $C_1,C_2$ be the two disks
  containing the edge $uv$. 
  Make a conforming split of $u$ into $u_1,u_2$ along $C_1$ and
  a conforming
  split of $v$ into $v_1,v_2$ along $C_2$ such that 
  both splits are conforming and $u_1$ and $v_1$ are
  adjacent in the resulting graph. 
  Now add the
  edge $u_2v_2$ to obtain \gp.
\item (split + cross) \label{split+cross} Let $u,v,w$ be vertices on a
  disk $C$ such that $u$ is not adjacent to $v$ or $w$. Perform a
  conforming split of $u$ into $u_1,u_2$, along $C$, with
  $u_1,u_2,v,w$ in that cyclic order on the new disk corresponding to
  $C$. Now add the edges $u_1v$ and $u_2w$ to obtain \gp.
\item (double split + cross) \label{dsplit+cross} Let $u,v$ be
  non-adjacent vertices on a disk $C$. Perform conforming splits of
  $u$ and $v$, into $u_1,u_2$ and $v_1,v_2$, respectively such that
  both splits are
  along $C$. Let $u_1,u_2,v_1,v_2$ appear in that cyclic order on the new
  disk corresponding to $C$. Now add the edges $u_1v_1$ and $u_2v_2$
  to obtain \gp.
\setcounter{opctr}{\value{enumi}}
\end{enumerate}
If \gp\ is obtained as in paragraph $i$ above, then we say that \gp\
is an \emph{\ext{$i$} of $G$ with respect to \D}. When the disk
system \D\ is implied by context, we may simply refer to an
\emph{\ext{$i$}} of $G$. We are now ready to state a preliminary form
of our main result, a counterpart of \thmref{thm:top-nonpl-ext},
with minors instead of topological minors.
A graph is a {\em prism} if it has exactly six vertices and its complement is
a cycle on six vertices.

\begin{theorem}
\label{thm:nonpl-ext}
Let $G$ be a weakly $4$-connected planar graph that is not a prism, let $H$ be a weakly
$4$-connected non-planar graph such that $G$ is isomorphic to a minor
of $H$, and let \D\ be the disk system in $G$ consisting of all
peripheral cycles. Then there exists an integer $i\in\{1,2,\ldots,7\}$
such that $H$ has a minor isomorphic to an \ext{$i$} of $G$ with
respect to \D.
\end{theorem}

\thmref{thm:top-nonpl-ext} is definitely easier to state than
\thmref{thm:nonpl-ext}. So what are the advantages of the latter
result? First, in the applications one is usually concerned with minors
rather than topological minors, and so \thmref{thm:nonpl-ext}
gives a more direct route to the desired results. Second, while the
number of types of outcome is larger in \thmref{thm:nonpl-ext},
in most cases the actual number of cases needed to examine will be
smaller. (Notice that, for instance, in \thmref{thm:top-nonpl-ext}
one must examine all $S$-paths between non-confluent ends, whereas in
\thmref{thm:nonpl-ext} one is only concerned with those between
non-confluent branch-vertices.)

Third, \thmref{thm:top-nonpl-ext} allows as an outcome an $S$-cross
on a cycle consisting of three segments. 
That is a drawback, which essentially means that in order for the theorem
to be useful the graph $G$ should have no triangles.
On the other hand, \thmref{thm:nonpl-ext} does not suffer from this
shortcoming and gives useful information even when $G$ has triangles.

Fourth, while a graph listed as an outcome of
\thmref{thm:top-nonpl-ext} may fail to be \wfc\ (and may do so in a
substantial way), an \ext{$i$} of a \wfc\ graph is again \wfc. This
has two advantages.  In the applications we are often seeking to prove
that \wfc\ graphs, with a minor isomorphic to some \wfc\ graph
embeddable in a surface $\Sigma$, that themselves do not embed into
$\Sigma$ have a minor isomorphic to a member of a specified list
$\mathcal{L}$ of graphs. In order to get a meaningful result we would
like each member of $\mathcal{L}$ to satisfy the same connectivity
requirement imposed on the input graphs.

From a more practical viewpoint, the advantage of maintaining the same
connectivity in the outcome graph is that the theorem can then be
applied repeatedly. That will become important when we consider a
generalization to arbitrary surfaces (that is, in the context of
theorems~\ref{thm:main} and \ref{thm:strong-main}). While a \wfc\
graph $G$ has at most one planar embedding, it may have several
embeddings in a non-planar surface $\Sigma$. Now one application of
the generalization of \thmref{thm:nonpl-ext} will dispose of one
embedding into $\Sigma$, but some other embedding might extend
naturally to those outcome graphs. So it may be necessary to apply the
theorem in turn to those outcome graphs in place of $G$. It will be
important that the outcomes of (the generalization of)
\thmref{thm:nonpl-ext} satisfy the same requirement as the input
graph. We can then apply such a theorem repeatedly till we get a list
of graphs that no longer embed in $\Sigma$ --- in other words, we
would have obtained the non-embeddable extensions of $G$. This will
be illustrated in \secref{sec:application}.

\section{Main Theorem}
\label{sec:main}


Our main theorem applies to arbitrary disks systems, at the expense of
having to add two outcomes. 
We also add a third additional outcome in order to allow $G$ to be a prism.
The extra outcomes are the following.
As before, let $G$ be a graph with
a cycle double cover \D, and let \gp\ be obtained by one of the operations
below. 
\begin{enumerate}
\setcounter{enumi}{\value{opctr}}
\item (non-separating triad) \label{triad} Let $x_1,x_2,x_3$ be three
  vertices of $G$ such that (i)~they are pairwise confluent, but not
  all contained in any single disk, and (ii)~$\{x_1,x_2,x_3\}$ is
  independent, and does not separate $G$. To obtain \gp, add a new
  vertex to $G$ adjacent to $x_1,x_2$ and $x_3$.
\item (non-conforming T-edge) \label{t-edge} Let a vertex $u$ and an
  edge $xy$ be such that (i)~$u$ is not confluent with the edge $xy$,
  but is confluent with both $x$ and $y$, (ii)~$u$ is not
  adjacent to either $x$ or $y$, and (iii)~$\{u,x,y\}$ does not
  separate $G$. Subdivide the edge $xy$ and join $u$ to the new
  vertex, to obtain \gp.
\item (enlargement of a prism) \label{prism} Let $G$ be a prism, and
  let \gp\ be obtained from $G$ by selecting two edges of $G$ that do
  not belong to a common peripheral cycle but both belong to a
  triangle, subdividing them, and joining the two new vertices by an
  edge.
\setcounter{opctr}{\value{enumi}}
\end{enumerate}
As before, if \gp\ is obtained as in paragraph $i$ above,
then we say that \gp\ is an
\emph{\ext{$i$} of   $G$ with respect to \D}.
Thus if $G$ is not a prism, then it has no \ext{\prism}, and
if $G$ is a prism, then its \ext{\prism} is unique, up to isomorphism.
The unique \ext{\prism} of the prism is known as $V_8$.

We also need to define an appropriate analogue of being non-planar in
the context of cycle double covers. That is the objective of this
paragraph and the next. Let $S$ be a subgraph of a graph $H$.  An
\emph{$S$-bridge of $H$} is a subgraph $B$ of $H$ such that either $B$
consists of a unique edge of $E(H)-E(S)$ and its ends, where the ends
belong to $S$, or $B$ consists of a component $J$ of $H\backslash
V(S)$ together with all edges from $V(J)$ to $V(S)$ and all their
ends.  For an $S$-bridge $B$, the vertices of $B\cap S$ are called the
\emph{attachments} of $B$. Let \D\ be a cycle double cover in $S$. We
say that \D\ is \emph{locally planar} in $H$ if the following
conditions are satisfied:
\begin{enumerate}
\renewcommand{\theenumi}{(\roman{enumi})}
\renewcommand{\labelenumi}{\theenumi}
\item for every $S$-bridge $B$ of $H$ there exists a 
    disk $C_B\in\D$ such that all the attachments of $B$ lie on
  $C_B$, and
\item for every disk $C\in\D$ the subgraph $\bigcup B\cup C$
  of $H$ has a planar drawing with $C$ bounding the unbounded
  face, where the big union is taken over all $S$-bridges $B$ of $H$ with
  $C_B=C$.
\end{enumerate}

Let $G$ have a weak disk system \D\ and $H$ have a minor isomorphic to
$G$. It is easy to see that there is an expansion $G'$ of $G$, such
that $G'$ is a topological minor of $H$.
We say that \D\ has a \emph{locally
  planar extension} into $H$ if:
\begin{enumerate}
\renewcommand{\theenumi}{(\roman{enumi})}
\renewcommand{\labelenumi}{\theenumi}
\item there exists a \emph{conforming} expansion $G'$ of $G$ such that
  a subdivision of $G'$ is a isomorphic to a subgraph $S$ of $H$, and
\item the weak disk system $\D'$ induced in $S$ by \D\ is locally planar
  in $H$.
\end{enumerate}


\noindent
We are now ready to state the main result.
\begin{theorem}
\label{thm:main}
Let $G$ and $H$ be \wfc\ graphs such that $H$ has a minor isomorphic
to $G$.
Let $G$ have a disk system \D\ that has no locally planar
extension into $H$. Then $H$ has a minor isomorphic to an \ext{$i$} of
$G$, for some $i\in\{1,\ldots,10\}$.
\end{theorem}
\bigskip

\noindent
Let us deduce \thmref{thm:nonpl-ext} from \thmref{thm:main}.
\bigskip

\noindent
{\bf Proof of \thmref{thm:nonpl-ext}, assuming \thmref{thm:main}.}
Let $G$  be as in \thmref{thm:nonpl-ext}, and let $i\in\{8,9,10\}$. 
By \thmref{thm:main} it suffices to show that
$G$ has no \ext{$i$} with respect to  the disk system consisting
of all peripheral cycles of $G$. 
This is clear when $i=10$, because $G$ is not a prism.
Thus we may assume for a contradiction
that  $i\in\{8,9\}$ and that such an \ext{$i$} exists. 
Let $u,x,y$ be the three vertices
of $G$ as in the definition of \ext{$i$}. Since every pair of vertices
among $u,x,y$ are confluent, it follows that $G\backslash \{u,x,y\}$
is disconnected, a contradiction.~\qed

\section{Outline of Proof}
\label{sec:outline}

\def\junk#1{}
\junk{
In fact, we prove a more general result, stated as
\thmref{thm:strong-main}, that applies to weak disk systems.  However,
this more general theorem requires three additional outcomes, and we
introduce them now. Let $J$ be a $K_4$ subgraph of a graph $G$, and
assume that each of the four peripheral cycles of $J$
are disks of $G$. In those circumstances we say that $J$ is a
\emph{detached $K_4$ subgraph} of $G$. Again, let $G$ be a graph with
a cycle double cover \D, and let \gp\ be obtained from $G$ by one of
the operations below.
\begin{enumerate}
\setcounter{enumi}{\value{opctr}}
\item (non-conforming double split) \label{ncdsplit} Let $v$ be a vertex
  of $G$ of degree at least 6. Partition the set 
  of edges incident with $v$ into three disjoint sets 
  $E_1,E_2,E_3$, each of size at least 2, such that 
  of the disks that use the vertex $v$, exactly two, say $D_1$ and $D_2$,
  use edges of two distinct sets $E_i$.
  Moreover, let $D_1$ and $D_2$ intersect precisely in the vertex $v$ and
  let them use no edge of $E_3$.
  Then \gp\ is obtained from $G$ by replacing $v$ with three vertices
  $v_1,v_2,v_3$, where $v_i$ is incident with edges in $E_i$ (for
  $i=1,2,3$) and $\gp[\{v_1,v_2,v_3\}]$ has edge-set $\{v_1v_3,v_2v_3\}$.
\item (simplex on detached $K_4$) \label{det-k4} 
  The graph \gp\ is obtained
  from $G$ by adding a new vertex adjacent to all four vertices of
  a detached $K_4$ subgraph of $G$.
\end{enumerate}
Similarly as before, if \gp\ is obtained as in paragraph $i$ above,
then we say that \gp\ is an \emph{\ext{$i$} of $G$ with respect to
  \D}. Thus if $G$ is not a prism, then it has no \ext{\prism}, and
if $G$ is a prism, then its \ext{\prism} is unique, up to isomorphism.
The unique \ext{\prism} of the prism is known as $V_8$.
}

The purpose of this section is to outline the proof of the main
theorem. 
Our main tool for the proof of \thmref{thm:main} will be its
counterpart for subdivisions, proved in \cite{RobSeyThoExt}. Before we
can state it we need one more definition.  Let $S$ be a subgraph of a
graph $H$, and let \D\ be a cycle double cover in $S$. Let $x\in
V(H)-V(S)$ and let $x_1,x_2,x_3$ be distinct vertices of $S$ such that
every two of them are confluent, but no disk of $S$ contains all
three. Let $L_1,L_2,L_3$ be three paths such that (i)~they share a
common end $x$, (ii)~they share no internal vertex among themselves or
with $S$, and (iii)~the other end of $L_i$ is $x_i$, for
$i=1,2,3$. The paths $L_1$, $L_2$, $L_3$ are then said to form an
\emph{$S$-triad}. The vertices $x_1,x_2,x_3$ are called the
\emph{feet} of the triad.
We are now ready to state our tool.
It is an immediate corollary of~\cite[Theorem~(4.6)]{RobSeyThoExt}.

\begin{theorem}[\cite{RobSeyThoExt}]
\label{thm:top-disk-ext}
Let $G$ be a graph with no vertices of degree two
that is not the complete graph on four vertices, let $H$ be a \wfc\
graph, let \D\ be a weak disk system in $G$, and let a subdivision of
$G$ be isomorphic to a subgraph of $H$. Then there exists a subgraph
$S$ of $H$ isomorphic to a subdivision of $G$ such that, letting $\D'$
denote the weak disk system induced in $S$ by \D, one of the following
conditions holds:
\begin{enumerate}
\item \label{top-jump} there exists an $S$-path in $H$ such that its
ends are not confluent in $S$, or
\item \label{top-cross} there exists a free $S$-cross in $H$ on some
  disk of $S$, or
\item \label{top-triad} the graph $H$ has an $S$-triad, or
\item the weak disk system $\D'$ is locally planar in $H$.
\end{enumerate}
\end{theorem}

Now let $G$, \D\ and $H$ be as in \thmref{thm:main}. It is easy to see
that there exists an expansion $G'$ of $G$ such that a subdivision of
$G'$ is isomorphic to a subgraph $S$ of $H$. (If $G$ itself is a
topological minor of $H$, then $G'=G$.) In \lemref{lem:nonconf-exp} we
prove that if $G'$ is a nonconforming expansion, then there exists a
\ext{\ncsplit} of $G$ that is isomorphic to a
minor of $H$. Thus from now on we may assume that $G'$ is a conforming
expansion of $G$. By \lemref{thm:top-disk-ext} applied to $S$ and $H$
we deduce that one of the outcomes of that lemma holds.  Notice that
those outcomes correspond to \ext{\jump}, \ext{\cross} and \ext{\triad},
respectively, except that in the enlargements the
vertices in question are required to be branch-vertices of $S$,
whereas in \lemref{thm:top-disk-ext} they are allowed to be interior
vertices of segments. We deal with this in \secref{sec:ext-of-G'} by
showing that each of the outcomes mentioned leads to a suitable
enlargement of $G'$. To be precise, at this point we settle for what
we call weak \triad- and weak \ext{\tedge}s, and in
\secref{sec:fromweak} show that these weak enlargements can be
replaced by ordinary enlargements, possibly of a different expansion
of $G$ and of a different kind. Finally, in \secref{sec:ext-of-G} we
complete the proof of \thmref{thm:main} by showing that the expansion
$G'$ can be chosen to be equal to $G$.

\section{Preliminaries}
\label{sec:prelim}

Let $G'$ be an expansion of a graph $G$. Then every vertex $v$ of
$G$ corresponds to a connected subgraph $T_v$ of $G'$. We call
$V(T_v)$ the \emph{branch-set corresponding to $v$}.

\begin{lemma}
\label{lem:tree}
Let $G'$ be an expansion of a graph $G$, let $u,v\in V(G)$ be distinct, 
and let $T_u,T_v$ be the corresponding subgraphs of $G'$.
Then $T_u$ and $T_v$ are induced subtrees of $G'$. If $u$ is adjacent
to $v$ then exactly one edge of $G'$ has one end in $V(T_u)$ and
the other in $V(T_v)$, and if  $u$ is not adjacent to $v$, then no
such edge exists.
\end{lemma}

An expansion of a \wfc\ graph may fail to be \wfc, but only in a limited
way. The next definition and lemma make that precise. 
Let $(A,B)$ be a nontrivial separation of order three in a graph $G$.
We say that $(A,B)$ is \emph{degenerate} if the vertices in $A\cap B$
can be numbered $v_1,v_2,v_3$ such that either
\begin{myenumerate}{}
\item $|A-B|=1$ and $A\cap B$ is an independent set, or
\item there exists a triangle $u_1u_2u_3$ in $G[A]$ such that
for $i=1,2,3$ the vertices $u_i$ and $v_i$ are either adjacent or
equal, $A\subseteq\{u_1,u_2,u_3,v_1,v_2,v_3\}$, and each edge of 
$G[A]$ is of the form $u_iv_i$ for $1\leq i\leq 3$ or $u_iu_j$ for
$1\leq i<j\leq 3$.
\end{myenumerate}
The following two lemmas are routine, and we omit the straightforward
proofs.

\begin{lemma}
\label{lem:3-cuts}
Let $G$ be an expansion of a \wfc\ graph. Then $G$ is
3-connected, and if it is not a prism, then for every nontrivial
separation $(A,B)$ of $G$ of order three, exactly one of  $(A,B)$,
 $(B,A)$ is degenerate.
\end{lemma}
\bigskip


\begin{lemma}
\label{lem:branch-sets-for-cut}
Let $G'$ be expansion of a \wfc\ graph $G$, let $(A,B)$ be a
degenerate separation of $G$ of order three satisfying condition
{\rm(2)} of the definition of degenerate separation, and let
$u_1,u_2,u_3,v_1,v_2,v_3$ be as in that condition. Then for at least
two integers $i\in \{1,2,3\}$ either $u_i=v_i$ or $u_iv_i$ is a new
edge of $G'$.
\end{lemma}

We now show that a non-conforming expansion of $G$ must have a minor
isomorphic to a \ext{\ncsplit} of $G$. 

\begin{lemma}
\label{lem:nonconf-exp}
Let \D\ be a disk system in a graph $G$, and 
let $G'$ be a non-conforming expansion of $G$. 
Then $G'$ has a minor isomorphic to a \ext{\ncsplit} of $G$.
\end{lemma}

\emph{Proof:}
We may assume that for every new edge $e$ of $G'$ the graph $G'/e$ is
a conforming expansion of $G$. We shall refer to this as the
\emph{minimality} of $G'$. We will prove that $G'$ is  a
\ext{\ncsplit} of $G$.

Let \ghat\ be an expansion of $G$ such that $G'$ is obtained from
\ghat\ by splitting a vertex $v$ into $v_1$ and $v_2$. By the
minimality of $G'$ this split is non-conforming, and \ghat\ is
a conforming expansion of $G$. If $G=\ghat$, then $G'$ is a \ext{\ncsplit}
of $G$, and so we may assume that $G\ne\ghat$. Let $e$ be a new
edge of \ghat. If $e$ is not incident with $v$, then $G'/e$ is
a non-conforming expansion of $\ghat/e$, contrary to the minimality
of $G'$. Now let us consider $e$ as an edge of $G'$. From the
symmetry between $v_1$ and $v_2$ we may assume that $e$ is incident
with $v_2$ in $G'$; let $v_3$ be its other end.
The split of the vertex $v$ of the graph \ghat\ into $v_1$ and $v_2$
violates \ref{split:1} or \ref{split:2}. But it does not violate
\ref{split:1}, for otherwise the same violation occurs in the analogous
split of $\ghat/e$, contrary to the minimality of $G'$.
Thus the split of the vertex $v$ of the graph \ghat\ into $v_1$ and $v_2$
satisfies \ref{split:1}; let $D_1$ and $D_2$ be the corresponding disks.
It follows that the  disks violate \ref{split:2}, but they do not
do so for the corresponding split in $\ghat/e$. It follows that 
$e\in E(D_1)\cap E(D_2)$.
Thus the split that creates \ghat\ from $\ghat/e$ is also along $D_1$
and $D_2$. Let $f$ be an edge incident with $v_2$ in $G'$ that is not
$e$ or the new edge $v_1v_2$ of $G'$.
It follows from~\ref{ax:2}
by considering the edge $f$ that either the split that creates 
 \ghat\ from $\ghat/e$ or the split that creates $G'/e$ from $\ghat/e$
is non-conforming, contrary to the minimality of $G'$.~\qed
\bigskip


The following lemma will be useful.

\begin{lemma}
\label{lem:two-disks}
Let $G'$ be a conforming expansion of a graph $G$ with respect to a
weak disk system \D, and let $\D'$ be the weak disk system induced in
$G'$ by \D. Let $qr$ be a new edge of $G'$, and let the vertex $p\in
V(G')-\{q,r\}$ share distinct disks $D_q,D_r$ of $G'$ with $q$ and
$r$, respectively, such that $D_r$ does not contain $q$. Then $p$ is
adjacent to $r$ and the disks $D_q,D_r$ both contain the edge $pr$.
\end{lemma}

\emph{Proof:}
The disks of $G/qr$ that correspond to $D_q$ and $D_r$ share $p$ and
the new vertex of $G/qr$, say $w$.  By (D3) $p$ is adjacent to 
$w$ in $G/qr$ and the edge $pw$ belongs to both those disks.
By \lemref{lem:tree} the vertex $p$ is adjacent to exactly one
of $q,r$.  But $q\notin V(D_r)$, and hence $p$ is adjacent to $r$
and $D_q,D_r$ both contain the edge $pr$, as desired.~\qed
\bigskip

We end this section with a lemma about fixing separations in
\wfc\ graphs, a special case of a lemma from \cite{JohThoSplit}.
First
some additional notation: when a
graph $G$ is a minor of a graph $H$, we say that an
\emph{embedding}
  $\eta$ of $G$ into $H$
is a mapping with domain $V(G) \cup E(G)$ as follows. $\eta$ maps
vertices $v \in G$ to connected subgraphs $\eta(v)$ of $H$, with
distinct vertices being mapped to disjoint vertex-disjoint subgraphs.
Further, $\eta$ maps edges $uv$ of $G$ to paths $\eta(uv)$ in $H$ with
one end in $\eta(u)$ and the other in $\eta(v)$, and otherwise
disjoint from $\eta(w)$ for any vertex $w$ of $G$. Also, for edges $e
\neq e'$ of $G$, if $\eta(e)$ and $\eta(e')$ share a vertex, then it
must be an end  of both the paths.
\begin{lemma}
\label{lem:ifc-ext}
Let $G_1$ be a graph isomorphic to a minor of a \wfc\ graph
$H$.
Let $P=\{p_1,p_2\}$, $Q=\{q_1,q_2,q_3\}$ and $R$ be such that
$(P,Q,R)$ is a partition of $V(G_1)$, and $G_1$ has all possible
edges between $P$ and $Q$, and no edge with both ends in $Q$.
Further, suppose $R$ has at least two vertices, and that $(P \cup Q,
Q \cup R)$ is a (non-trivial) 3-separation of $G_1$.
Then $H$ has a minor isomorphic to a graph $G_1^+$ that is obtained
from $G_1$ by
\begin{enumerate}
\item adding an edge between $p_i$ and $r$ for some $i \in \{1,2\}$
  and $r \in R$, or
\item splitting $q_j$ for some $j \in \{1,2,3\}$ into vertices $q_j^1$
  and $q_j^2$ such that $q_j^1$ is adjacent to $p_1$ and $q_j^2$ is
  adjacent to $p_2$
\end{enumerate}
\end{lemma}

\emph{Proof:}
Call an embedding $\eta$ of $G_1$ into $H$ \emph{minimal} if for every
embedding $\eta'$ of $G$ into $H$,
\[\sum_{j=1}^3 |E(\eta(q_j))| \leq \sum_{j=1}^3 |E(\eta'(q_j))|\]
In particular, if $\eta$ is minimal, $\eta(q_j)$ is
a tree for every $j$. Further, we say that a vertex $q_j$ is
\emph{good} for $\eta$ if the paths $\eta(p_1q_j)$ and $\eta(p_2q_j)$
are vertex-disjoint (in other words, their ends in $\eta(q_j)$ are
distinct).

Consider a minimal embedding $\eta$ of $G_1$ into $H$.  Suppose there
exists a $q_j$ that is good for $\eta$. For $i=1,2$, let $p_i'$ be the
endpoint of $\eta(p_iq_j)$ in $\eta(q_j)$. Let $e$ be an edge in the
unique path between $p_1'$ and $p_2'$ in $\eta(q_j)$, and let
$T_1,T_2$ be the two subtrees obtained by deleting $e$ from
$\eta(q_j)$, such that $p_i' \in T_i$ for $i=1,2$. For $i=1,2$, define
$N_i$ to be the set of neighbors $r \in R$ of $q_j$ in $G$ that
$\eta(rq_j)$ has an endpoint in $T_i$. Now $N_1,N_2$ are non-empty by
the minimality of $\eta$. (If, say, $N_1$ were empty, then we could
replace $\eta(q_j)$ by $T_2$ and modify $\eta(p_1q_j)$ accordingly to
get a better embedding $\eta'$, a contradiction.) It is easy to see
that conclusion~2 of the lemma is satisfied, with the neighborhoods of
$q_j^1$ and $q_j^2$ being $N_1 \cup \{p_1\}$ and $N_2 \cup \{p_2\}$,
respectively.

Hence we may assume that there is no minimal embedding of $G$ into $H$
with a vertex in $Q$ being good for it. Let $\eta$ be an embedding of
$G$ into $H$. For $j=1,2,3$, there exist vertices $t_j$ such that both
$\eta(p_1q_j)$ and $\eta(p_2q_j)$ have $t_j$ as an end. Define $J_1$
as the union of $\eta(p_i)$, $i=1,2$ and of $\eta(e)$ for all edges
$e$ with at least one end in $P$.  Define $J_2$ as the union of
$\eta(v)$ for $v \in Q \cup R$ and of $\eta(e)$ for every edge $e$ of
$G$ with both ends in $Q \cup R$. Now $V(J_1) \cap V(J_2) =
\{t_1,t_2,t_3\}$. Since $H$ is \wfc, there is a path in $H$ with ends
$a \in V(J_1) \setminus V(J_2)$ and $b \in V(J_2) \setminus V(J_1)$,
and otherwise disjoint from $J_1 \cup J_2$. If $b$ belongs to
$\eta(q_j) \setminus t_j$ for some $j$, then we can modify $\eta$ to
get a minimal embedding where $q_j$ is a good vertex, which is a
contradiction. Thus $b$ belongs to $\eta(r)$ for some $r \in R$ or $b$
is an internal vertex of $\eta(e)$ for an edge $e$ of $G$ that has an
end in $R$ (recall that $Q$ is an independent set). In either case, it
is easy to see that conclusion~1 holds.\qed

\section{The Enlargements of an Expansion of \texorpdfstring{$G$}{G}}
\label{sec:ext-of-G'}

Let $G$ and $H$ be as in \thmref{thm:main}.  In order to apply
\thmref{thm:top-disk-ext} we select an expansion $G'$ of $G$ such
that a subdivision of $G'$ is isomorphic to a subgraph of $H$. By
\lemref{lem:nonconf-exp} we may assume that $G'$ is a conforming
expansion. In this section we prove three lemmas, one corresponding to
each of the first three outcomes of \thmref{thm:top-disk-ext}.
The lemmas together almost imply that the conclusion of
\thmref{thm:main} holds for $G'$. The reason for the word almost
is that for convenience we allow a weaker form of \ext{\triad}s and
\ext{\tedge}s.

The weaker form of \ext{\tedge}s is defined as follows.
Let $G$ be a graph with a cycle double cover \D, and let $u,x,y\in V(G)$,
where $x$ and $y$ are adjacent and $u$ is not confluent with the edge
$xy$. Let \gp\ be obtained from $G$ by subdividing the edge $xy$ and
adding an edge joining the new vertex to $u$. We say that \gp\ 
is a weak \ext{\tedge} of $G$. 
Later, in \lemref{lem:weak-t-edge},
we show how to move from a weak \ext{\tedge} to a \ext{\tedge} or another
useful outcome. Our first lemma deals with the first outcome of
\thmref{thm:top-disk-ext}.

\begin{lemma}
\label{lem:ext-of-G'no1}
Let $G,H$ be graphs such that $G$ is connected, has at least five
vertices and no vertices of degree two.  Let \D\ be a weak disk system
in $G$, let $S$ be a subgraph of $H$ isomorphic to a subdivision of
$G$, and let $P$ be an $S$-path in $H$ such that its ends are not
confluent in the weak disk system \D$'$ induced in $S$ by \D. Then
$H$ has a minor isomorphic to a \ext{\jump}, \ext{\ncsplit} or a weak
\ext{\tedge} of $G$.
\end{lemma}

\emph{Proof:} Let $s,t$ be the ends of $P$.  If both $s$ and $t$ are
branch-vertices in $S$, then $H$ has a minor isomorphic to a
\ext{\jump} of $G$, and we are done.  If one of $s$ and $t$ is a
branch-vertex and the other is an internal vertex of a segment of $S$,
then $H$ has a minor isomorphic to a weak \ext{\tedge} of $G$, as
desired.

Thus we may assume that $s$ and $t$ are internal vertices of two
different segments $Q_1$ and $Q_2$ of $S$, respectively. Let $Q_1$
correspond to an edge $uv\in E(G)$, and let $Q_2$ correspond to an
edge $xy\in E(G)$.
Now, if $u$ is not confluent with the edge $xy$, then $H$ has a minor
isomorphic to a weak \ext{\tedge} of $G$, and the lemma holds. Thus,
we may assume that $u$ shares a disk $D_1$ with the edge $xy$. By
symmetry, we get a disk $D_2$ shared by $v$ and the edge $xy$, and
disks $D_3,D_4$ that the edge $uv$ shares with vertices $x$ and $y$
respectively. (The disks $D_i$ may not be pairwise distinct.)

The disks $D_1$ and $D_3$, however, must be distinct, since the
vertices $s,t$ are not confluent. Notice, however, that they share the
vertices $u$ and $x$.  It follows that $u,v,x,y$ are pairwise
distinct, for if $v=y$, say, then $u, v=y, x$ all belong to $V(D_1\cap
D_3)$, and hence $D_1=D_3$ by \ref{ax:3}, a contradiction. By
\ref{ax:3} this implies that $u$ is adjacent to $x$ in $G$ and the
intersection of $D_1$ and $D_3$ is precisely the edge $ux$. In other
words, the vertices $u$ and $x$ must be adjacent in $G$, and $D_1,D_3$
are precisely the two disks containing the edge $ux$. By a similar
argument, it follows that $u$ and $y$ are adjacent, and $D_1,D_4$ are
precisely the two disks containing the edge $uy$. Thus $u$ is adjacent
to each of $v,x,y$ in $G$, and the edges $uv,ux,uy$ are pairwise
confluent.

By symmetry, we get similar conclusions about the vertices
$v,x,y$. Thus $G[u,v,x,y]$ is a detached $K_4$ subgraph of $G$.
Since $G$ has at least five vertices and is connected, we may assume,
without loss of generality, that $u$ has a neighbor in $G$ outside of
the set $\{v,x,y\}$. Let $N$ be the set of all such neighbors of $u$.
But then delete the edges of the segment of $S$ corresponding to the
edge $ux$ and contract the edges of the subpath of $Q_2$ between $t$
and the end corresponding to $x$. It follows that $H$ has a minor
isomorphic to a graph obtained from $G$ by splitting $u$ corresponding
to the partition $\{\{v,x\},N\cup\{y\}\}$ of its neighbors. This split
is non-conforming since the disks $D_1$ and $D_4$ violate
condition~\ref{split:2} in the definition of a conforming split. Hence
$H$ has a minor isomorphic to a \ext{\ncsplit} of $G$.~\qed


\begin{lemma}
\label{lem:ext-of-G'no2}
Let $G,H$ be graphs such that $H$ is \wfc, and $G$ is connected, has
at least 5 vertices and has no vertices of degree two. Let \D\ be a
weak disk system in $G$, let $S$ be a subgraph of $H$ isomorphic to a
subdivision of $G$, such that $\D$ induces the weak disk system $\D'$
in $S$.  Further, let there exist a free $S$-cross on some disk of
$S$.  Then $H$ has a minor isomorphic to a \ext{\cross} or a
\ext{\ncsplit} or a weak \ext{\tedge} of $G$.
\end{lemma}

\emph{Proof:}
Let the free cross consist of paths $P_1,P_2$, in a disk
$C'$ of $S$, that corresponds to a disk $C$ of $G$.  We
shall call the paths $P_1,P_2$ the \emph{legs} of the cross. Recall that
the ends of $P_1,P_2$ are called the \emph{feet} of the cross.

If $C$ has at least four vertices, then we claim that $H$ has a minor
isomorphic to a \ext{\cross} of $G$.
We define an auxiliary bipartite graph $B$, with the vertex set being
the set of feet of the cross and the set of branch-vertices of $S$
that belong to $C'$. A foot $f$ and a branch-vertex $b$ are adjacent
if one of the subpaths of $C'$ with ends $f$ and $b$ includes no feet
or branch-vertices in its interior. Since the cross is free, it
follows from Hall's bipartite matching theorem that $B$ has a complete
matching from the set of feet to the set of branch vertices (in other
words, one that matches each of the feet). By contracting the edges of
the paths that correspond to this matching, we deduce that $H$ has a
minor isomorphic to a \ext{\cross} of $G$, as desired.

Hence we may assume that $C$ is in fact a triangle on vertices
$u_1,u_2$ and $u_3$, say. For $i=1,2,3$, if $u_i$ has degree 3 in $G$,
then define $v_i$ to be its third neighbor (that is, the neighbor not
in $C$). Otherwise, define $v_i=u_i$. Let
$u'_1,u'_2,u'_3,v'_1,v'_2,v'_3$ be the corresponding vertices of
$S$. Let $Q_i$ denote the segment of $S$ corresponding to the edge
$u_iv_i$ if $u_i\ne v_i$ and let $Q_i$ be the null graph otherwise,
and let $A=V(C'\cup P_1\cup P_2)$ and $B=(V(S)-V(C'\cup Q_1\cup
Q_2\cup Q_3))\cup\{v'_1,v'_2,v'_3\}$.  There exist three
vertex-disjoint paths in $H$ linking $\{u'_1,u'_2,u'_3\}$ to
$\{v'_1,v'_2,v'_3\}$. Since $H$ is \wfc, it follows that there is no
3-cut in $H$ separating $A$ from $B$.  Hence, by a variant of Menger's
theorem, $H$ contains four vertex-disjoint paths $L_1,\ldots L_4$
linking $\{v'_1,v'_2,v'_3,y\}$ to $\{u'_1,u'_2,u'_3,x\}$ (not
necessarily in that order), where $x \in A$ and $y \in B$. We assume
the numbering of the paths is such that for $i=1,2$ and $3$, $L_i$ has
end $v'_i \in B$. (The remaining path $L_4$ then has end $y \in B$.)
We wish to define a suitable vertex $w\in V(G)$.
If $y$ is a branch-vertex, then let $w$ be the corresponding vertex of $G$;
otherwise $y$ is an internal vertex of a segment of $H$.
By Lemmas~\ref{lem:tree} and~\ref{lem:branch-sets-for-cut} at least
one end of that segment, say $w'$, does not belong to $\{v_1,v_2,v_3\}$,
and we let $w$ be the vertex of $G$ that corresponds to $w'$.


We may assume that $x \in V(C')$. If not, then we may contract edges
suitably in $P_1$ or $P_2$ such that the vertex corresponding to $x$,
after the contraction, lies on $C'$. (Note that this contraction does
not affect the graph $S$, neither does it destroy the
cross.)

Relabel the vertices $u'_1,u'_2,u'_3,x$ as $a,b,c,d$, in the order in
which they appear on $C'$ (in some orientation), such that $L_4$ joins
$d$ to $y$. (Note that $d$ need not be the same as $x$.)
%
Let $(d,a,b)$ denote the interior vertices of the subpath of $C'$
with ends $d$ and $b$ that includes $a$ in its interior, and let
$(d,c,b)$ be defined analogously.

We claim that there is a leg of the cross with feet $f,g$ such that
$f \in (d,a,b)$ and $g \in (d,c,b)$.
Since the cross is free, there exists a leg with foot in $(d,a,b)$.
We may assume the other foot of this leg does not belong to
$(d,c,b)$, but then the other leg of the cross satisfies the claim. 

Choose a leg as above such that there is no foot between $f$ and $a$,
and no foot between $g$ and $c$. (Such a choice must be possible, due
to the freeness of the cross.) Let the other leg of the cross have
feet $h,i$, such that $b$ and $h$ are joined by a subpath of the cycle
$C'$ that is disjoint from $\{f,g\}$. By contracting disjoint
subpaths of $C'$ with ends $(a,f)$, $(c,g)$, and $(b,h)$ respectively,
%
it follows that $H$ has a minor isomorphic to the graph
$G'$ obtained from $G$ by adding a new vertex $z$ adjacent to $u_1$,
$u_2$, $u_3$ and $w$.

If $w$ is not confluent with the edge $u_1u_2$ then $G'\setminus
u_1u_2\setminus zu_3$ is isomorphic to a weak \ext{\tedge} of $G$, and
we are done. Thus we may assume that $w$ is confluent with the edge
$u_1u_2$, and by symmetry, with the edges $u_2u_3$ and $u_1u_3$ as
well.
%
%
It follows similarly as in the proof of \lemref{lem:ext-of-G'no1}
that $G[u_1,u_2,u_3,w]$ is a detached $K_4$ subgraph of $G$. Since $G$
is connected, and $|V(G)|\geq 5$, we may assume, without loss of
generality, that $u_1$ has a neighbor in $G$ outside that set. It
follows that a graph obtained from $G$ by a non-conforming split of
$u_1$ is isomorphic to a minor of $H$.\qed

\medskip
We now define the weaker form of \ext{\triad}s.
Let $G$ be a graph with a cycle double cover \D, and let $x_1,x_2,x_3$
be vertices of $G$ such that no disks contains all three.
Let \gp\ be obtained from $G$ by adding a vertex with neighborhood
$\{x_1,x_2,x_3\}$. We say that \gp\
is a weak \ext{\triad} of $G$.
Our third lemma deals with the third outcome of
\thmref{thm:top-disk-ext}.


\begin{lemma}
\label{lem:ext-of-G'no3}
Let $G,H$ be graphs such that $G$ is connected, has at least five
vertices and no vertices of degree two. Let \D\ be a weak disk system
in $G$, let $S$ be a subgraph of $H$ isomorphic to a subdivision of
$G$, and let there exist an $S$-triad in $H$. Then $H$ has a minor
isomorphic to an \ext{$i$} of $G$ for $i= \jump\mbox{ or }\ncsplit$,
or a weak \ext{$i$} of $G$ for $i= \triad\mbox{ or }\tedge$.
\end{lemma}

\emph{Proof:} We proceed by induction of $|E(H)|$. Let the $S$-triad
be $L_1,L_2,L_3$, and let its feet be $x_1,x_2,x_3$. If each $x_i$ is
a branch-vertex of $S$, then $S\cup L_1\cup L_2\cup L_3$ gives rise to
a minor of $H$ isomorphic to a weak \ext{\triad}, as desired. We may
therefore assume that $x_3$ is an internal vertex of a segment $Q_3$
of $S$ with ends $u_3$ and $v_3$. Let $f$ be an edge of $Q_3$. By
induction applied to $G$, $H/f$, and $S/f$, we may assume that $f$ is
incident with $x_3$ and one end of $Q_3$, say $u_3$, and that there
exists a disk $D_1$ in $G$ containing $x_1,x_2,u_3$. Similarly, we may
assume that there exists a disk $D_2$ in $G$ containing
$x_1,x_2,v_3$. Then $D_1\ne D_2$, because otherwise $D_1=D_2$ includes
the segment $Q_3$ by \ref{ax:3}, and hence each of $x_1,x_2,x_3$, a
contradiction. Since $x_1$ and $x_2$ belong to $D_1\cap D_2$, it
follows from \ref{ax:3} that $x_1$ and $x_2$ belong to a common
segment $Q$ of $S$.

Let $S'$ be obtained from $S$ by replacing $Q[x_1,x_2]$ by $L_1\cup
L_2$. Applying \lemref{lem:ext-of-G'no1} to $G$, $H$, $S'$ and the
$S'$-path $L_3$, the lemma now follows.~\qed

\bigskip
Using \thmref{thm:top-disk-ext} we can
summarize Lemmas~\ref{lem:ext-of-G'no1}--\ref{lem:ext-of-G'no3} as follows.

\begin{lemma}
\label{lem:ext-of-G'}
Let $G,H$ be \wfc\ graphs, let $G$ have a disk system \D\
with no locally planar extension into $H$,
and let 
a subdivision of $G$ be isomorphic to a subgraph of $H$.
Then $H$ has a minor isomorphic to
\begin{myenumerate}[\roman]{}
\item an \ext{$i$} of $G$ for some $i \in \{\jump, \cross, \ncsplit\}$, or
\item a weak \ext{$i$} of $G$ for some $i \in \{\triad, \tedge\}$.
\end{myenumerate}
\end{lemma}

\emph{Proof:} 
Let $G,H$ and \D\ be as stated. By \thmref{thm:top-disk-ext}
we deduce that there exists a subgraph $S$ of $H$
isomorphic to a subdivision of $G$ such that the induced disk system
 in $S$ satisfies one of the outcomes of 
\thmref{thm:top-disk-ext}. But then (i) or (ii)
of this lemma holds by Lemmas~\ref{lem:ext-of-G'no1}--\ref{lem:ext-of-G'no3}.~\qed

\section{From Weak Enlargements to Enlargements}
\label{sec:fromweak}

The purpose of this section is to replace  weak enlargements by
enlargements in \lemref{lem:ext-of-G'}(ii). We start with a special
case of weak \ext{\tedge}s.

\begin{lemma}
\label{lem:weak-t-edge-prelim}
Let $G$ be a graph with a cycle double cover \D, let
\gp\ be a weak \ext{\tedge} of $G$, and let $u,x,y$ be as in
the definition of weak \ext{\tedge}. If $G\backslash\{u,x,y\}$ is
connected,  then
\gp\ has a minor isomorphic to an \ext{$i$} of $G$ for some $i \in
\{\jump, \ncsplit, 
 \tedge\}$.
\end{lemma}

\emph{Proof:}
Let $z$ be the new vertex of \gp\ that resulted from the subdivision of
the edge $xy$.
If $u$ and $x$ are not confluent, then contracting the edge $xz$ of
\gp\ produces a \ext{$1$} of $G$, and so the lemma holds. Thus
we may assume that  $u$ and $x$ are confluent, and, by symmetry,
we may assume that $u$ and $y$ are confluent.
If $u$ is not adjacent to $x$ or $y$, then \gp\ is a \ext{$9$} of $G$,
and the lemma holds. Thus, from the symmetry, we may assume that
$u$ is adjacent to $x$.
The edges $xy$ and $xu$ are not confluent, for otherwise $u$ is
confluent with the edge $xy$, contrary to what a weak \ext{\tedge}
stipulates. But then deleting the edge $xu$ from \gp\ yields a graph
isomorphic to a \ext{\ncsplit} of $G$ --- more specifically, a graph
obtained by a non-conforming split of the vertex $x$.~\qed

\begin{lemma}
\label{lem:vertexadd}
Let $G$ be an expansion of a \wfc\ graph, 
let $H$ be a \wfc\ graph, let \D\ be a weak disk system 
in $G$, let $v$ be a vertex of $G$ of degree three and let
$u,x,y$ be the neighbors of $v$. Let \gp\ be the graph obtained from
$G$ by adding a new vertex $z$ adjacent to $u,x,y$ and deleting all
edges with both ends in $\{u,x,y\}$. If $H$ has a minor isomorphic to
\gp, then $H$ has a minor isomorphic to an \ext{$i$} of $G$ for some
$i\in\{\jump, \cross, \ncsplit, \splitjump,\tedge,\prism\}$.
\end{lemma}

\emph{Proof:}
Since $G$ is an expansion of a \wfc\ graph, \lemref{lem:3-cuts}
implies that at most one edge of $G$ has both ends in $\{u,x,y\}$.
Thus we may assume that $u$ is not adjacent to $x$ or $y$.
Since $v$ has degree three, it follows from \ref{ax:1} and \ref{ax:3} that
if $x$ is adjacent to $y$, then
the triangle $vxy$ is a disk in $G$. 
We can apply \lemref{lem:ifc-ext} to
$\gp=G_1$ and $H$, with $P=\{v,z\}$, $Q=\{u,x,y\}$ and
$R=V(\gp)-(P\cup Q)$. 
  From
the lemma, using the symmetry between $x$ and $y$, and the symmetry
among $x$, $y$ and $u$ if $x$ is not adjacent to $y$, 
we get the following three cases:

\textbf{Case 1:} $H$ has a minor isomorphic to a graph \gpp\ that is
obtained from \gp\ by adding an edge between a vertex $p \in P$ and a
vertex $r \in R$. Note that the vertices $v$ and $z$ are symmetric for
the application of \lemref{lem:ifc-ext}. Hence we may assume that
$p=v$. Now if $r$ is not confluent with $v$ in $G$, then \gpp\ above
has a minor isomorphic to a \ext{\jump} of $G$. Thus we may assume
that $r$ is confluent with $v$ in $G$. Furthermore, we may assume,
without loss of generality, that the disk $D_3$ shared by $r$ and $v$
contains the edges $vu$ and $vy$. (Note that $v$ has degree 3 in $G$.)
On the disk $D_3$, the vertices $u,v,y$ and $r$ occur in that cyclic
order. Now in \gpp, contracting the edge $yz$ gives a cross in the disk
$D_3$ with arms $uy$ and $rv$. In other words, \gpp, and hence $H$,
has a minor isomorphic to a \ext{\cross} of $G$, as desired.

\textbf{Case 2:} 
The vertices $x$ and $y$ are adjacent in $G$ and
$H$ has a minor isomorphic to a graph \gpp\ that is
obtained from \gp\ by splitting the vertex $x$
into $x_1$ and $x_2$, with $x_1$ adjacent to $v$ and $x_2$ adjacent to $z$.
Let $N_i$ be the neighbors of $x_i$ in \gpp\ other than $v,z,x_1,x_2$.
The neighborhood of $x$ in $G$ is thus $N_1 \cup N_2 \cup \{v,y\}$. In
$G$, let $D_4$ be the disk that contains the edge $xy$, other than the
triangle $vxy$. The disk $D_4$ must contain a vertex in either $N_1$ or
$N_2$, and from the symmetry between $v$ and $z$ we may assume that
it contains a vertex in $N_1$. Then, in \gpp, delete
the edge $uz$ and contract the edge $x_2z$. This gives a graph that is
a \ext{\ncsplit} of $G$ (non-conforming split of $x$, with the disks
$vxy$ and $D_4$ violating condition~\ref{split:2} in the definition of
a conforming split), as desired.

\textbf{Case 3:}  $H$ has a minor isomorphic to a graph \gpp\ that is
obtained from \gp\ by splitting the vertex $u$ into $u_1$ and
$u_2$, with $u_1$ adjacent to $v$ and $u_2$ adjacent to $z$.
Let $N_i$ be the set of neighbors of $u_i$ other than $v,z,u_1,u_2$.
Thus in $G$, the neighborhood of $u$ is $N_1 \cup N_2
\cup \{v\}$.

Let $D_1$ be the disk in $G$ shared by the edges $xv$ and $vu$, and
$D_2$ be the disk in $G$ shared by the edges $yv$ and $vu$. The disks
$D_1$ and $D_2$ both contain exactly one vertex each from $N_1\cup
N_2$. Let us assume first that $|N_2|\ge2$.
Contract the edge $xz$ in \gpp, and if $x$ is not adjacent to $y$ in
$G$, then delete also the resulting edge $xy$ to obtain a graph $G_1$,
and let $G_2$ be the graph obtained from $G_1$ by further deleting the
edge $u_2x$. Now $G_2$ is isomorphic to a graph obtained from $G$ by
splitting the vertex $u$ into $u_1$ and $u_2$. If this split is
non-conforming, then $G_2$ is a \ext{\ncsplit} of $G$, and we are
done. Otherwise, the split is not along $D_1$ or $D_2$, and from the
symmetry we may assume it is not along $D_1$. Thus $G_1$ is a
\ext{\splitjump} of $G$. (Note that in $G$, $u$ and $x$ are
non-adjacent, and hence non-consecutive on $D_1$.)  This completes the
case when $|N_2|\ge2$.

From the symmetry we may therefore assume that $|N_1|=|N_2|=1$. Thus
the degree of $u$ in $G$ is three.  For $i=1,2$ let $N_i=\{n_i\}$.  We
may assume that the edge $un_i$ belongs to the disk $D_i$. It follows
that the vertex $x$ and edge $un_2$ are not confluent in $G$, for if
some disk $D$ contained both of them, then the intersection $D\cap
D_1$ would violate \ref{ax:3}, because $u$ is not adjacent to $x$.
The graph $G_1$ from the previous paragraph is a weak \ext{\tedge} of
$G$, and so by \lemref{lem:weak-t-edge-prelim} we may assume that
$G\backslash\{x,u,n_2\}$ is disconnected. Since $u$ has degree three,
the \wfcity\ of $G$ implies that $n_1$ has degree three and its
neighbors are $x,u,n_2$. Since $G\backslash \{n_2,y\}$ is connected,
we deduce that $G$ is isomorphic to the prism, and \gpp\ is isomorphic
to a \ext{\prism} of $G$, as desired.~\qed \bigskip

Now we are ready to eliminate weak \ext{\tedge}s.

\begin{lemma}
\label{lem:weak-t-edge}
Let $G$ be an expansion of a \wfc\ graph, let \D\ be a weak disk system
in $G$, and let \gp\ be a weak \ext{\tedge} of $G$ such that \gp\ is
isomorphic to a minor of a \wfc\ graph $H$. Then $H$ 
has a minor isomorphic to an \ext{$i$} of $G$ for some $i \in
\{\jump, \cross, \ncsplit, \splitjump, \tedge,\prism\}$.
\end{lemma}

\emph{Proof:}
Let $u,x,y$ be as in the definition of weak \ext{\tedge}.
By \lemref{lem:weak-t-edge-prelim} we may assume
that $G\backslash\{u,x,y\}$ is disconnected. Since 
$x$ is adjacent to $y$ and $G$ is an expansion of a \wfc\ graph, 
\lemref{lem:3-cuts}
implies that the neighborhood of some vertex $v$ of $G$ is precisely the
set $\{u,x,y\}$. Thus \gp\ is as described in \lemref{lem:vertexadd},
and the conclusion follows from that lemma.~\qed
\bigskip

We now turn to weak \ext{\triad}s. In order to save effort we prove
a weaker analogue of \lemref{lem:weak-t-edge}, the following.

\begin{lemma}
\label{lem:weak-triad}
Let $G_1$ be an expansion of a \wfc\ graph $G$, let \D\ be a weak disk
system in $G$, and let $\gp$ be a weak \ext{\triad} of $G_1$ such that
$\gp$ is isomorphic to a minor of a \wfc\ graph $H$. Then there exists
an expansion $G_2$ of $G$ obtained from $G_1$ by contracting a
possibly empty set of new edges such that $H$ has a minor isomorphic
to an \ext{$i$} of $G_2$ for some $i \in \{\jump, \cross, \ncsplit,
\splitjump, \triad, \tedge,\prism\}$.
\end{lemma}
 
\emph{Proof:} We proceed by induction on $|E(G_1)|$. Let \gp\ be
obtained from $G_1$ by adding a vertex joined to $v_1,v_2,v_3$. If
some edge of $G_1$ has both ends in the set $\{v_1,v_2,v_3\}$, then by
deleting that edge we obtain a graph isomorphic to a weak \ext{\tedge}
of $G_1$, and the lemma follows from \lemref{lem:weak-t-edge}. Thus
we may assume that $\{v_1,v_2,v_3\}$ is an independent set in
$G_1$. We may also assume that every pair of vertices in
$\{v_1,v_2,v_3\}$ is confluent, for otherwise \gp\ has a minor
isomorphic to a \ext{\jump} of $G$, and the lemma holds. Thus we may
assume that $G\backslash\{v_1,v_2,v_3\}$ is disconnected, for
otherwise $G^+$ is an \ext{\triad} of $G_1$.

Let $(A,B)$ be a non-trivial separation of $G$ with $A\cap
B=\{v_1,v_2,v_3\}$. By \lemref{lem:3-cuts} we may assume that
$(A,B)$ is degenerate. If $|A-B|=1$, then the lemma follows from
\lemref{lem:vertexadd}. Thus we may assume that $|A-B|\geq 2$. Let
$v_1,v_2,v_3,u_1,u_2,u_3$ be as in the definition of degenerate. Since
$\{v_1,v_2,v_3\}$ is independent, we may assume from the symmetry that
$u_1\neq v_1$ and $u_2\neq v_2$. Now one of $u_1v_1,u_2v_2$ is a new
edge of $G_1$, and so we may assume the former is. Thus $\gp/u_1v_1$
is a weak \ext{\triad} of $G_1/u_1v_1$, and hence the lemma follows by
the induction hypothesis applied to the graph $G_1/u_1v_1$.~\qed
\bigskip

The lemmas of this section allow us to upgrade \lemref{lem:ext-of-G'}
to the following.

\begin{lemma}
\label{lem:ext-of-G'-new}
Let $G,H$ be \wfc\ graphs, let $G$ have a disk system \D\
with no locally planar extension into $H$,
and let $G'$ be a conforming expansion of $G$ such that
a subdivision of $G'$ is isomorphic to a subgraph of $H$.
Then there exists a conforming expansion $G''$ of $G$ obtained from
$G'$ by contracting a possibly empty set of 
new edges such that, letting $\D''$ denote
the weak disk system induced in $G''$ by \D, the graph
$H$ has a minor isomorphic to
an \ext{$i$} of $G''$ with
respect to $\mathcal{D}''$ for some
$i \in \{\jump, \cross, \ncsplit, \splitjump, \triad, \tedge, \prism\}$.
\end{lemma}

\emph{Proof:}
By \lemref{lem:ext-of-G'} we may assume that a weak \ext{\triad}
or a weak \ext{\tedge} of $G'$ is isomorphic to a minor of $H$.
By Lemmas~\ref{lem:weak-t-edge} and~\ref{lem:weak-triad} there exists
a required conforming expansion $G''$ of $G$ such that 
$H$ has a minor isomorphic to an \ext{$i$} of $G''$ for some $i \in
\{\jump, \cross, \ncsplit, \splitjump, \triad, \tedge, \prism\}$.~\qed

\section{Proof of the Main Theorem}
\label{sec:ext-of-G}

\lemref{lem:ext-of-G'-new} gives an \ext{$i$} of an expansion 
$G''$ of $G$.  Our final objective is to show that we can choose
$G''=G$.  We break the proof into several lemmas depending on the value
of $i$.

\begin{lemma}
\label{lem:op-jump}
Let $G$ and $H$ be \wfc\ graphs,
and let \D\ be a weak disk system in $G$ with no locally
planar extension into $H$. Let $G'$ be a
conforming expansion of $G$ such that $H$ has a minor isomorphic to a
\ext{\jump} of $G'$. Then $H$ has a minor isomorphic to an \ext{$i$} of
$G$ for some $i\in\{\jump,\ncsplit,\splitjump,\dsplitjump\}$.
\end{lemma}

\emph{Proof:} We may assume that $G'$ is as stated in the lemma, and
subject to that, it is minor-minimal. By hypothesis, $H$ has a minor
isomorphic to \gp, a graph obtained from $G'$ by adding an edge
between two vertices $x$ and $y$ that are not confluent. Let $e$ be a new
edge of $G'$.  By the minimality of $G'$, it
follows that
\begin{enumerate}
\renewcommand{\theenumi}{(\roman{enumi})}
\renewcommand{\labelenumi}{\theenumi}
\item one end of $e$ must be in $\{x,y\}$, and 
\item the other end of $e$ must be confluent with the vertex in
  $\{x,y\}$ other than the one above.
\end{enumerate}
Recall that branch-sets of an expansion were defined at the beginning of
\secref{sec:prelim}.
Thus all branch sets that are disjoint from $\{x,y\}$ are singleton
sets. Let $T_p$ and $T_q$ be the branch sets 
corresponding to vertices $p,q\in V(G)$ such 
that they contain $x$ and $y$
respectively ($p$ and $q$ may be identical). We claim that the degree
of $x$ in the branch set containing it is at most one (that is, $x$ is
a leaf of the tree $G'[T_p]$). Suppose not; hence $x$ has (at least)
two neighbors $x_1$ and $x_2$ in $T_p$. By (ii) above, $y$ shares
disks $D_1$ and $D_2$ of $G'$ with $x_1$ and $x_2$ respectively. 
Then $x\notin V(D_1\cup D_2)$, for $x,y$ are not confluent.  It
follows that $D_1\ne D_2$, for otherwise $D_1$ is not a cycle in
$G'/x_1x/x_2x$, and yet $D_1$ corresponds to a disk in $G$.
Also, $y$ is not adjacent to both $x_1$ and $x_2$, by
\lemref{lem:tree}.  But then contracting edges $xx_1$ and
$xx_2$ violates Axiom~\ref{ax:3} in $G$. This proves the claim. Thus $x$, and
by symmetry $y$, are leaf vertices in $G'[T_p]$ and $G'[T_q]$
respectively.

If $p=q$, then it follows that $T_p=T_q$ must be a path of length 2,
with a middle vertex $z$. Let $D'_1,\D'_2$ be the two disks in $G'$
that include the edge $xz$, and let ${D'_3},{D'_4}$
be the two disks that include the edge $yz$. Note that, since $x$ and
$y$ are not confluent in $G'$, all four disks above are distinct. Let
$D_1,D_2,D_3,D_4$ be the corresponding disks in $G$. Let $N_1,N_2$ be
the partition of the set of neighbors of $p$ in $G$, corresponding to
the partition $\{x,y\},\{z\}$ of $V(T_p)$. Clearly, $N_1$ has at least
two vertices, but so does $N_2$, by Axiom~\ref{ax:3} applied to
$\Dat,\Dbt$. In \gp\ (which has the edge $xy$), contract the edge
$xy$. This gives a graph \gpp\ that can be obtained from $G$ by
splitting $p$ with respect to the partition $N_1,N_2$ of its
neighbors. This split is non-conforming, since the disks
$D_1,\ldots,D_4$ violate condition~\ref{split:1} in the definition of
a conforming split.  Thus \gpp\ is a \ext{\ncsplit} of $G$, as desired.

If $p \neq q$, then from (i) and (ii) above, $T_p$ is either $\{x\}$
or $\{x,x_1\}$. By symmetry, $T_q$ is either $\{y\}$ or $\{y,y_1\}$.
If $T_p$ and $T_q$ are both singletons, then clearly $G' = G$ and we
are done.

Suppose exactly one of the two branch sets, say $T_q$, is a singleton,
and $T_p$ consists of $\{x,x_1\}$, where $x_1$ shares a disk $D$ with
$y$ in $G'$. If $x_1$ and $y$ are not adjacent, then \gp\ is a
\ext{\splitjump} of $G$, and we are done. Thus we may assume that
$x_1$ and $y$ are adjacent, and hence by Axiom~\ref{ax:3}, they are
consecutive in $D$. Let $D_1,D_2$ be the two disks in $G'$ containing
the edge $xx_1$. They are both distinct from $D$, since $x$ and $y$
are not confluent in $G'$. By Axiom~\ref{ax:3} applied to $D_1$ and $D_2$,
the vertex $x_1$ has at least two neighbors in $G'$ other $x$ and $y$. Now in
\gp\ (which contains the edge $xy$), delete the edge $x_1y$. This
gives a graph $\widetilde G$ obtained from $G$ by splitting $p$ in the same way
as in $G'$, except that $y$ is adjacent to $x$ rather than
$x_1$. Further, it is a non-conforming split, as the disks $D$, $D_1$
and $D_2$ violate condition~\ref{split:1} in the definition of a
conforming split. Thus $\widetilde G$, which is isomorphic to a minor of $H$, is
a \ext{\ncsplit} of $G$, and we are done.

Finally, suppose $T_p = \{x,x_1\}$ and $T_q = \{y,y_1\}$, where $x$
shares a disk $D'_1$ with $y_1$ and $y$ shares a disk $D'_1$ with
$x_1$. Let $D_1,D_2$ be the corresponding disks in $G$. Since $x$ and
$y$ are not confluent in $G'$, $D'_1$ does not contain $y$ and $D'_2$
does not contain $x$. (In particular, $D'_1$ and $D'_2$ are distinct.)
Apply \lemref{lem:two-disks} to $\ghat=G'/xx_1$, with the
vertices $p,y,y_1$ in that graph corresponding to $p,q,r$ in the
lemma. Thus the (conforming) split of the vertex $q$ in $G$ that
produces \ghat\ is along $D_2$, and $D_2$ is one of the disks
containing the edge $pq$ in $G$. Also, since $x$ and $y$ are not
confluent in $G'$, the (conforming) split of $p$ in \ghat\ that
produces $G'$ must be along $D_1$, and $D_1$ is the other disk in $G$
containing $pq$. It now follows that \gp\ is a \ext{\dsplitjump} of
$G$.
This finishes the proof of the lemma.\qed

\begin{lemma}
\label{lem:op-cross}
Let $G$ and $H$ be \wfc\ graphs, and let \D\ be a weak disk system in
$G$ with no locally planar extension into $H$. Let $G'$ be a
conforming expansion of $G$ such that $H$ has a minor isomorphic to a
\ext{\cross} of $G'$. If $G'\ne G$, then there exists a conforming
expansion $G''$ of $G$ obtained from $G'$ by contracting at least one
new edge such that $H$ has a minor isomorphic to an \ext{$i$} or a
weak \ext{\tedge} of $G''$ for some
$i\in\{\cross,\splitcross,\dsplitcross\}$.
\end{lemma}

\emph{Proof:} We may assume that $G'$ is as stated in the lemma, and
subject to that, it is minor-minimal. By hypothesis, there are
vertices $u,v,x,y$ appearing on a disk $C'$ in $G'$, in that cyclic
order, such that $H$ has a minor isomorphic to a graph  obtained
from $G'$ by adding the edges $ux$ and $vy$. Let $C$ be the cycle in
$G$ corresponding to $C'$. 
The minimality of $G'$ implies that every new edge of $G'$ has both
ends in $\{u,v,x,y\}$, and hence it belongs to $C'$ by \ref{ax:3}.
We may therefore assume that $uv$ is a new edge of $G'$.  We claim that if
$v$ is adjacent to $x$, then the lemma holds.  To prove this claim
suppose that $v$ and $x$ are adjacent in $G'$, and let $G_1=G^+\backslash
vx$.  If $v$ has degree three in $G'$, then $G_1$ is isomorphic to a weak
\ext{\tedge} of $G'/uv$ (the new edge is $yv$; notice that $y$ is
not confluent with the edge of $G'/uv$ that is being subdivided
by \ref{ax:3}), and hence the lemma holds.  Thus we may assume that
$v$ has degree at least four in $G'$.  In that case $G_1$ is isomorphic to
a \ext{\splitjump} of $G'/uv$, for a graph isomorphic to $G_1$ can be
obtained by a conforming split of the new vertex of $G'/uv$, not
along $C'$, and joining one of the new vertices to $y$.  This proves
our claim, and hence we may assume that $v$ is not adjacent to $x$.
By symmetry we may also assume that $u$ is not adjacent to $y$.

If $uv$ is the only new edge of $G'$, then $G'$ is a \ext{\splitcross}
of $G$, and the lemma holds.  Thus we may assume that $G'$ has another
new edge, and so that edge must be $xy$ and there are no other new edges.
It follows that $G'$ is a \ext{\dsplitcross} of $G$, and so the
lemma holds.~\qed
\bigskip

\begin{lemma}
\label{lem:op-tedge}
Let $G$ and $H$ be graphs, let \D\ be a weak disk system in $G$,
and let $G'$ be a
conforming expansion of $G$ such that $H$ has a minor isomorphic to a 
\ext{\tedge} \gp\ of $G'$. 
If $G'\ne G$, then there exists a conforming expansion $G''$ obtained from
$G'$ by contracting at least one new edge such that
$H$ has a minor isomorphic to a \ext{\ncsplit} or a weak 
\ext{\tedge} of $G''$.
\end{lemma}

\emph{Proof:}
Let $u,x,y\in V(G')$ be such that
\gp\ is obtained from $G'$ by subdividing the edge $xy$ and joining
the new vertex to $u$, and let $f$ be a new edge of $G'$.
Then $f\ne xy$, for otherwise \lemref{lem:two-disks}
implies that $u$ is confluent with the edge $xy$, a
contradiction. We may assume that $f$ is
incident with $u$, and that contracting $f$ makes the new vertex
 confluent with the edge $xy$, for otherwise $\gp/f$ is a weak
\ext{\tedge} of $G'/f$, and the lemma holds.
Hence the other end $v$ of $f$ must share a disk
$D_{1}$ with the edge
$xy$. Since $u$ is not confluent
with $xy$, $D_{1}$ does not contain $u$. Let $D_{2}$ and $D_{3}$ be
disks shared by $u$ and $x$, and by $u$ and $y$, respectively. These
three disks are pairwise distinct, since $u$ is not confluent with the
edge $xy$ in $G'$. Now apply \lemref{lem:two-disks} with $x$ as the
vertex $p$, and $u,v$ as the vertices $q,r$ respectively. It follows
that $v$ and $x$ are adjacent, and that $D_{1}$ and $D_{2}$ are the
two disks containing the edge $vx$. Apply \lemref{lem:two-disks}
again, this time with $y$ in place of $x$. It follows that the edges
$vu,vx$ and $vy$ are covered twice each by the three disks
$D_{1},D_{2}$ and $D_{3}$. In particular, $D_1$ is a triangle.

If $f'\ne f$ is a new edge of $G'$, then by what we have shown about
$f$ it follows that $f'$ is incident with $u$ and its other end belongs
to a disk $D_1'$ that contains the edge $xy$. Since $D_1$ is a triangle
consisting of $x,y$ and an end of $f$, we see that $D_1'\ne D_1$.
But the disks that correspond to $D_1$ and $D_1'$ in $G'/f/f'$ have three
vertices in common, contrary to \ref{ax:3}. Thus $f$ is the only new edge
of $G'$, and hence $G=G'/f$. Let $p$ be the new vertex of $G=G'/f$.

Since \gp\ is a \ext{\tedge} of $G'$, the graph $G'\backslash\{u,x,y\}$
is connected, and hence
 $v$ has a neighbor outside $\{u,x,y\}$. (In
fact, it must then have at least three neighbors outside $\{u,x,y\}$.)
Let $z$ be the new vertex of \gp\ created by subdividing the edge
$xy$. The graph $\gp\backslash vx/xz$ is isomorphic to a graph obtained from
$G$ by splitting $p$ into two vertices.
This split is non-conforming, since the two disks in
$G$ that contain $py$ violate condition~\ref{split:2} in the
definition of a conforming split. Thus $H$ has a minor isomorphic to a
\ext{\ncsplit} of $G$.  This finishes the proof of the lemma.~\qed

\begin{lemma}
\label{lem:op-triad}
Let $G$ and $H$ be graphs, let \D\ be a disk system in $G$,
and let $G'$ be a
conforming expansion of $G$ such that $H$ has a minor isomorphic to an
\ext{\triad} \gp\ of $G'$.
If $G'\ne G$, then there exists a conforming expansion $G''$ obtained from
$G'$ by contracting at least one new edge such that
$H$ has a minor isomorphic to a \ext{\ncsplit}
or a weak \ext{\triad}  of $G''$.
\end{lemma}

\emph{Proof:}
Let \gp\ be obtained from $G'$ by adding a vertex adjacent to $x_1,x_2,x_3$,
and let $f$ be a new edge of $G'$.
We may assume that upon contracting $f$ the vertices that correspond to
$x_1,x_2,x_3$ belong to a common disk, for otherwise $\gp/f$ is a weak
\ext{\triad} of $G'/f$, and the lemma holds.
Thus $f$ is incident with at least one of 
$x_1,x_2,x_3$, say $x_1$, and there exists a disk \D\ in $G'$ that includes
$y,x_2,x_3$, where $y$ is the other end of $f$.

Apply \lemref{lem:two-disks} twice, once with $x_2$ as the vertex $p$,
and next with $x_3$ as the vertex $p$. In both applications, let $x_1$
and $y$ be the vertices $q$ and $r$ respectively. It follows that $y$
is adjacent to $x_2$ and $x_3$, and that $yx_1\in E(D_2\cap D_3)$,
$yx_2\in E(D\cap D_3)$ and $yx_3\in E(D\cap D_2)$.  Since \gp\ is a
\ext{\triad} of $G'$ the graph $G'\backslash\{x_1,x_2,x_3\}$ is
connected, and hence $y$ has degree at least four.  Let $N$ be the
neighbors of $y$ in $G'$ other than $x_1,x_2,x_3$.  Let $G'$ be
obtained from $G$ by splitting $x_1$ in such a way that the
neighborhood of one of the new vertices is $N$. Then $G'$ is
isomorphic to a minor of \gp, and it is a \ext{\ncsplit} of $G'$.
Thus the lemma follows from \lemref{lem:nonconf-exp}.~\qed


\bigskip

We are finally ready to state and prove 
\thmref{thm:main}, which we restate.

\begin{theorem}
\label{thm:strong-main}
Let $G$ and $H$ be \wfc\ graphs such that $H$ has a minor isomorphic
to $G$.  Let $G$ have a disk system \D\ that has no locally
planar extension into $H$. Then $H$ has a minor isomorphic to an
\ext{$i$} of $G$, for some $i\in\{1,2,\ldots,10\}$.
\end{theorem}

\emph{Proof:}
There exists an expansion of $G$ whose subdivision is isomorphic to
a subgraph of $H$. If this expansion is not conforming, then the
theorem holds by \lemref{lem:nonconf-exp}, and so we may assume
that the expansion is conforming.
By \lemref{lem:ext-of-G'-new} there exists a conforming expansion
$G'$ of $G$ such that $H$ has a minor isomorphic to an
\ext{$i$} \gp\ of $G'$ for some 
$i\in\{\jump,\cross,\ncsplit,\splitjump, \triad, \tedge,\prism\}$.
We may choose $G'$ and \gp\ such that $|E(G')|$ is minimum. 
If $i\in\{\jump,\splitjump\}$,
then \gp\ is isomorphic to a \ext{\jump} of a conforming expansion
of $G'$, and the theorem holds by \lemref{lem:op-jump}.
If $i=\ncsplit$, then the theorem holds by 
\lemref{lem:nonconf-exp}.
If $i=\prism$, then the minimality of $G'$ implies that $G=G'$,
and if $i\in\{\cross,\triad,\tedge\}$, then
the same conclusion follows from Lemmas~\ref{lem:op-cross},
\ref{lem:op-triad} and~\ref{lem:op-tedge}, respectively,
using Lemmas~\ref{lem:weak-t-edge} and~\ref{lem:weak-triad}.
Thus the theorem holds.~\qed


\section{An Application}
\label{sec:application}

In this section, we illustrate an application of
\thmref{thm:main}. Archdeacon \cite{arch-phd,arch-proj} proved that a
graph $H$ does not embed in the projective plane if and only if it has
a minor isomorphic to some graph in an explicitly constructed list of
35~graphs. One might hope that if we assume that $H$ is sufficiently
connected, then the list may be shortened. Mohar and Thomas (work in
progress) developed a strategy for a proof, but it will be a lengthy
project with several intermediate steps. Here we complete one such
step: under the assumptions that $H$ is \wfc\ and has a minor
isomorphic to the Petersen graph, \thmref{thm:petersen} below gives a
list of eight forbidden minors, each of which are \wfc.

\begin{figure}
\begin{tabular*}{\textwidth}{c@{\extracolsep{\fill}}c}
\includegraphics[height=0.19\textheight]{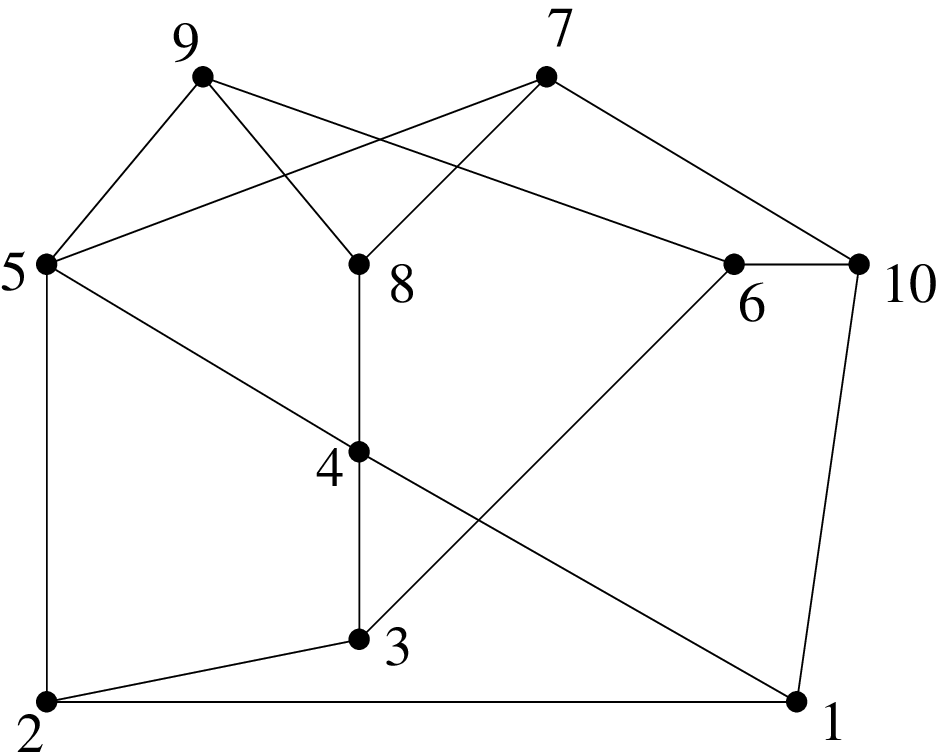}&%
\includegraphics[height=0.19\textheight]{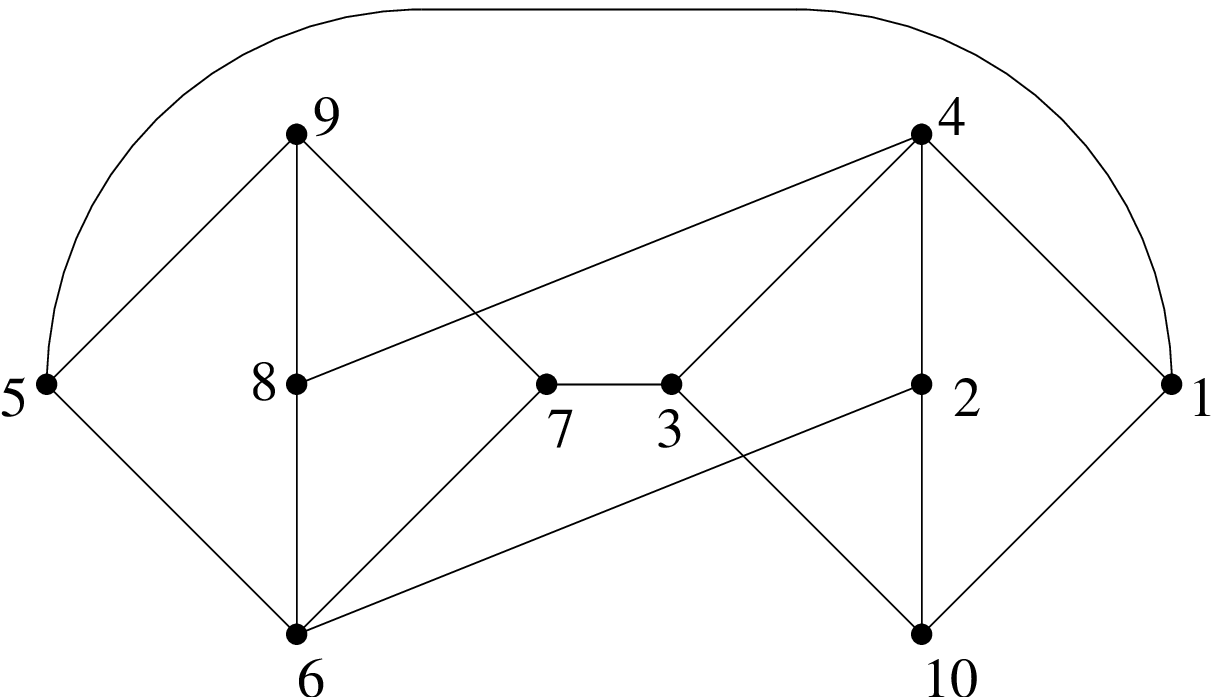}\\
\fsplit&\ffour\\
 & \\
\includegraphics[height=0.19\textheight]{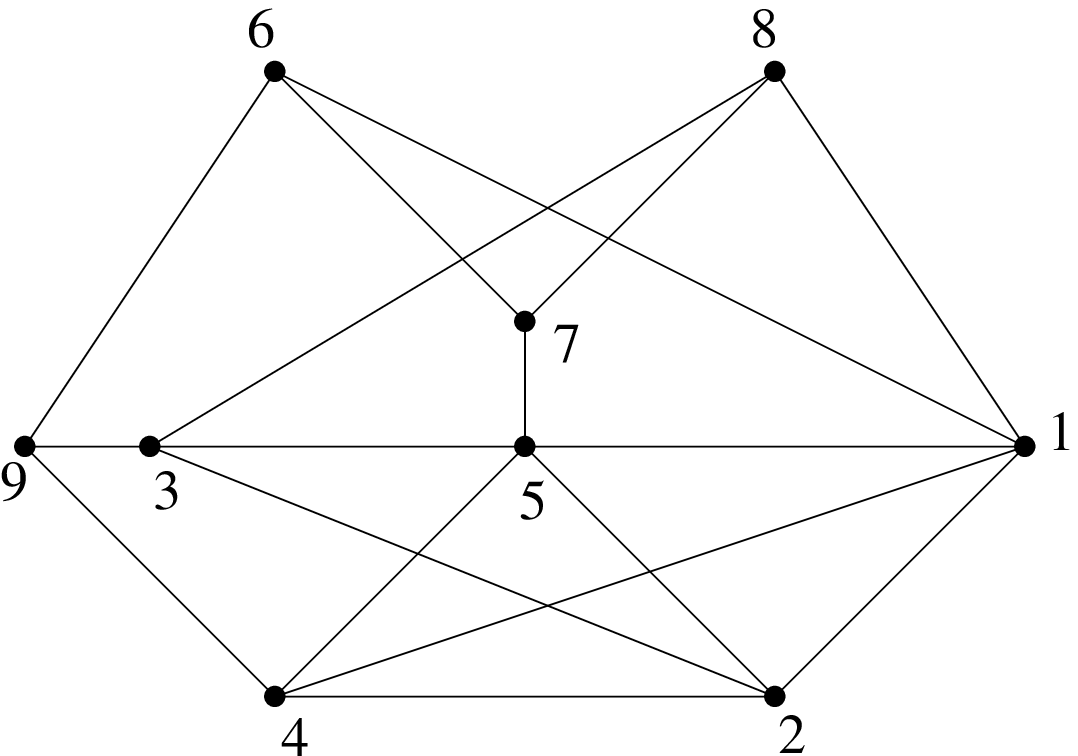}&%
\includegraphics[height=0.19\textheight]{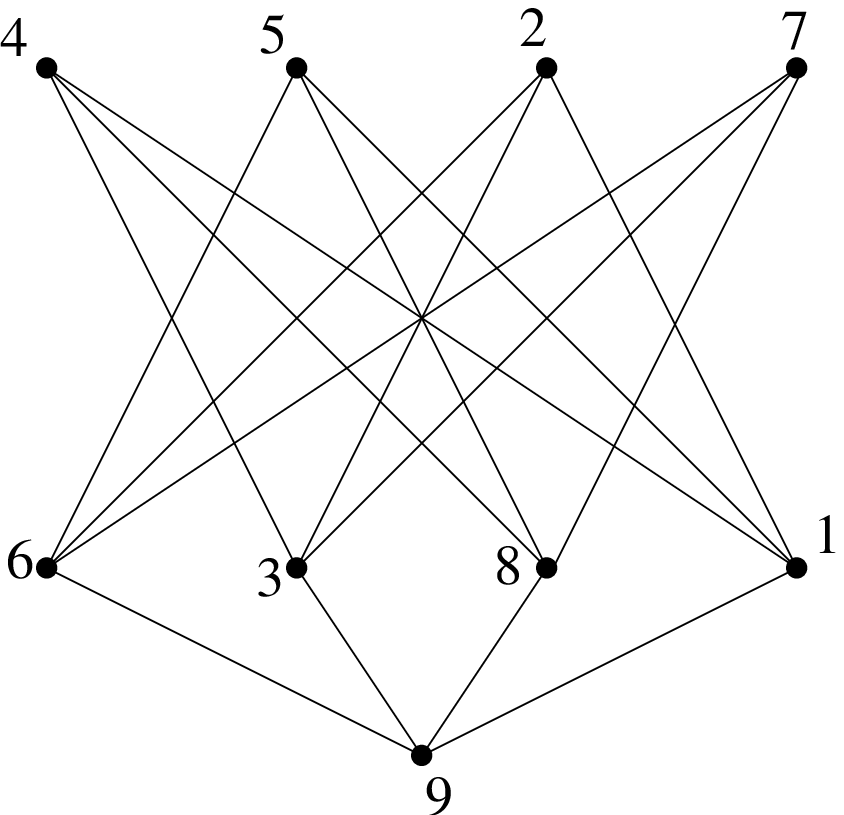}\\
\dsplit&\kffmfour\\
 & \\
\includegraphics[height=0.19\textheight]{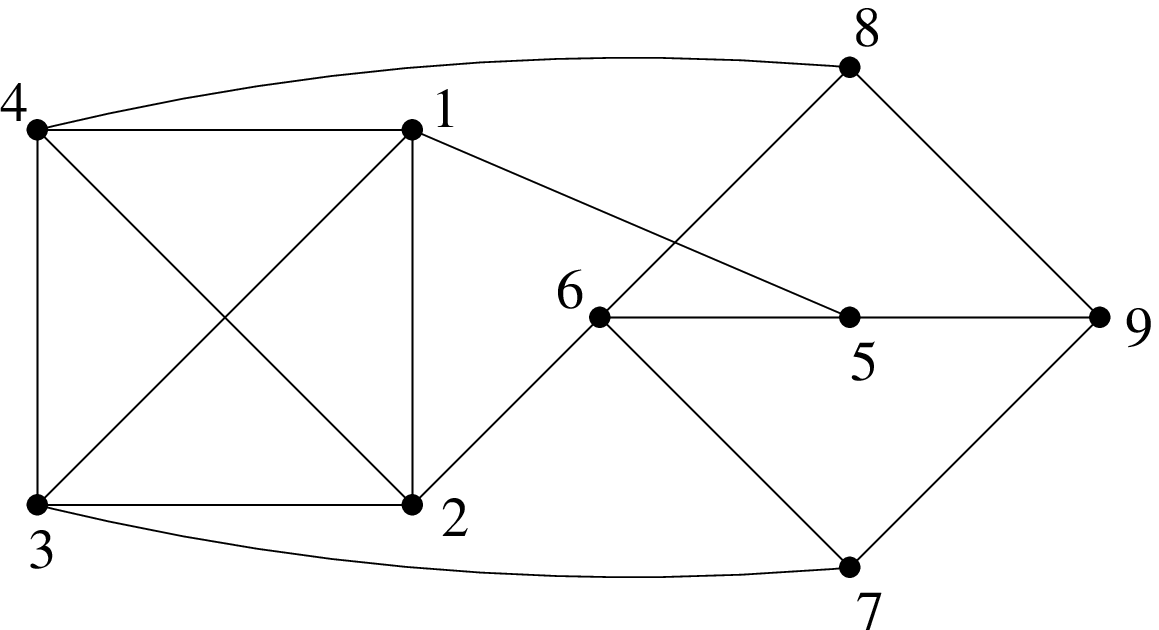}&%
\includegraphics[height=0.19\textheight]{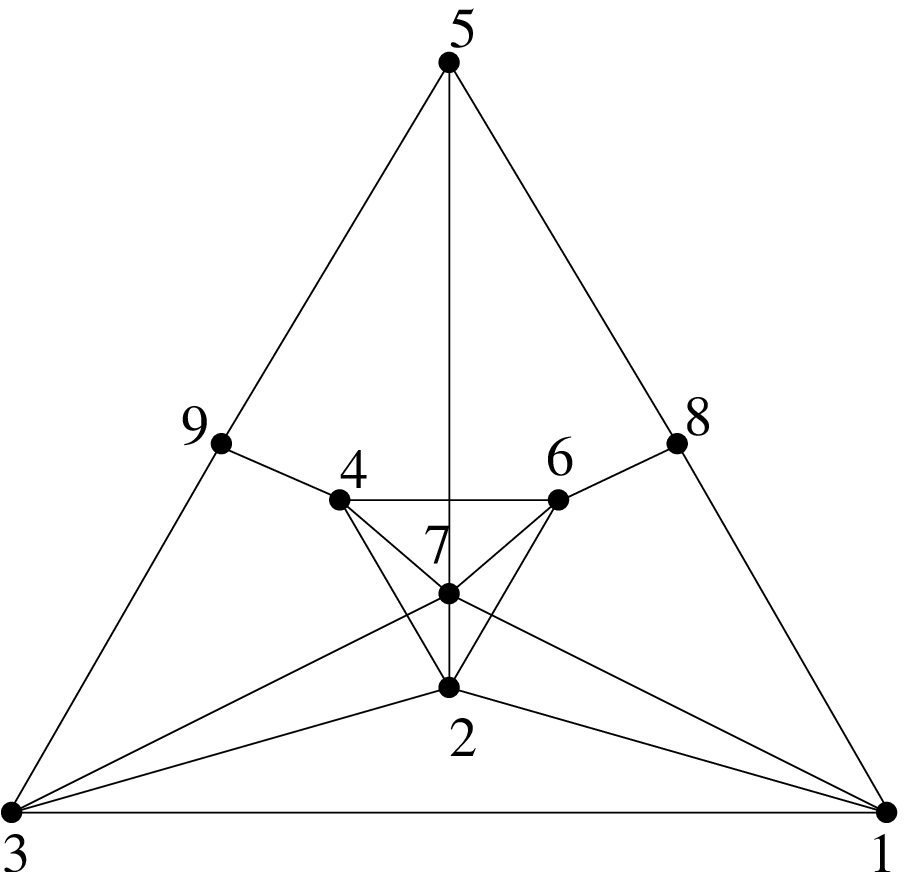}\\
\etwenty&\cthree\\
 & \\
\includegraphics[height=0.19\textheight]{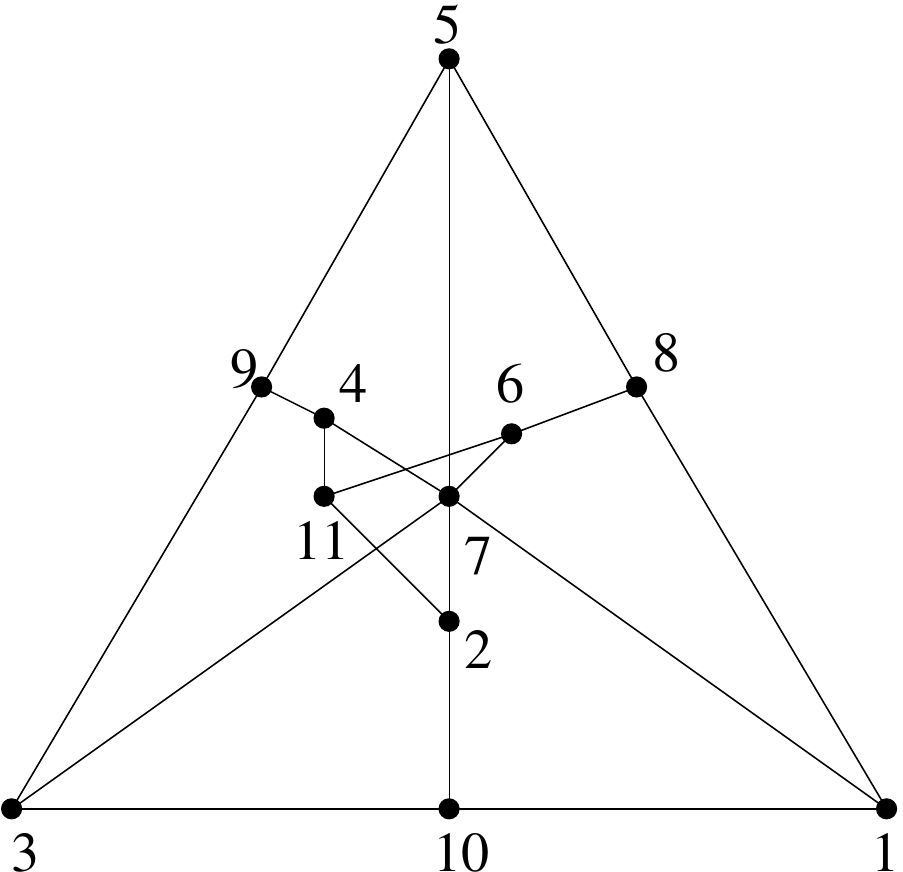}&%
\includegraphics[height=0.19\textheight]{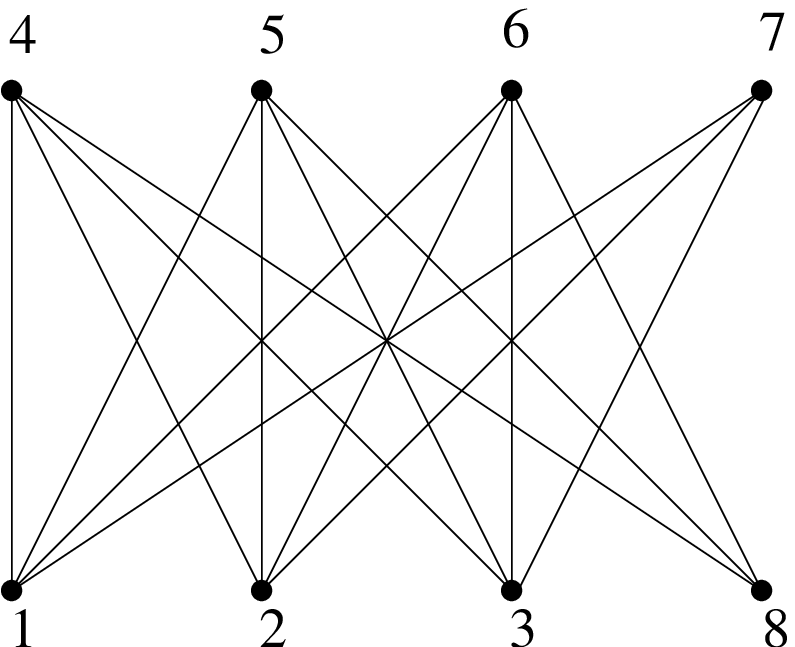}\\
\etwo&\kffminus
\end{tabular*}
\caption{The eight graphs of \thmref{thm:petersen}}
\label{fig:proj-list}
\end{figure}

\figref{fig:proj-list} shows these eight graphs (with a
vertex-labeling for each of them).
All of these graphs, with the exception of \fsplit\ and \dsplit,
appear in the list of 35 forbidden minors for the projective
plane. \fsplit\ and \dsplit, however, are obtained from two graphs in
that list ($F_1$ and $D_3$, respectively) by splitting exactly one
vertex. (The reason we list $\fsplit,\dsplit$ instead of $F_1,D_3$ is
that the latter two graphs are not \wfc.)

\begin{theorem}
\label{thm:petersen}
Let $H$ be a \wfc\ graph that has a minor isomorphic to the Petersen
graph. Then $H$ does not embed in the projective plane if and only if
it has a minor isomorphic to one of the eight graphs \fsplit, \ffour,
\dsplit, \kffmfour, \etwenty, \cthree, \etwo, or \kffminus\ shown in
\figref{fig:proj-list}.
\end{theorem}

Before we derive \thmref{thm:petersen} from \thmref{thm:main}, we
describe some notation that will be convenient in the proof.

\begin{figure}[ht]
  \centering
  \includegraphics[width=0.4\textwidth]{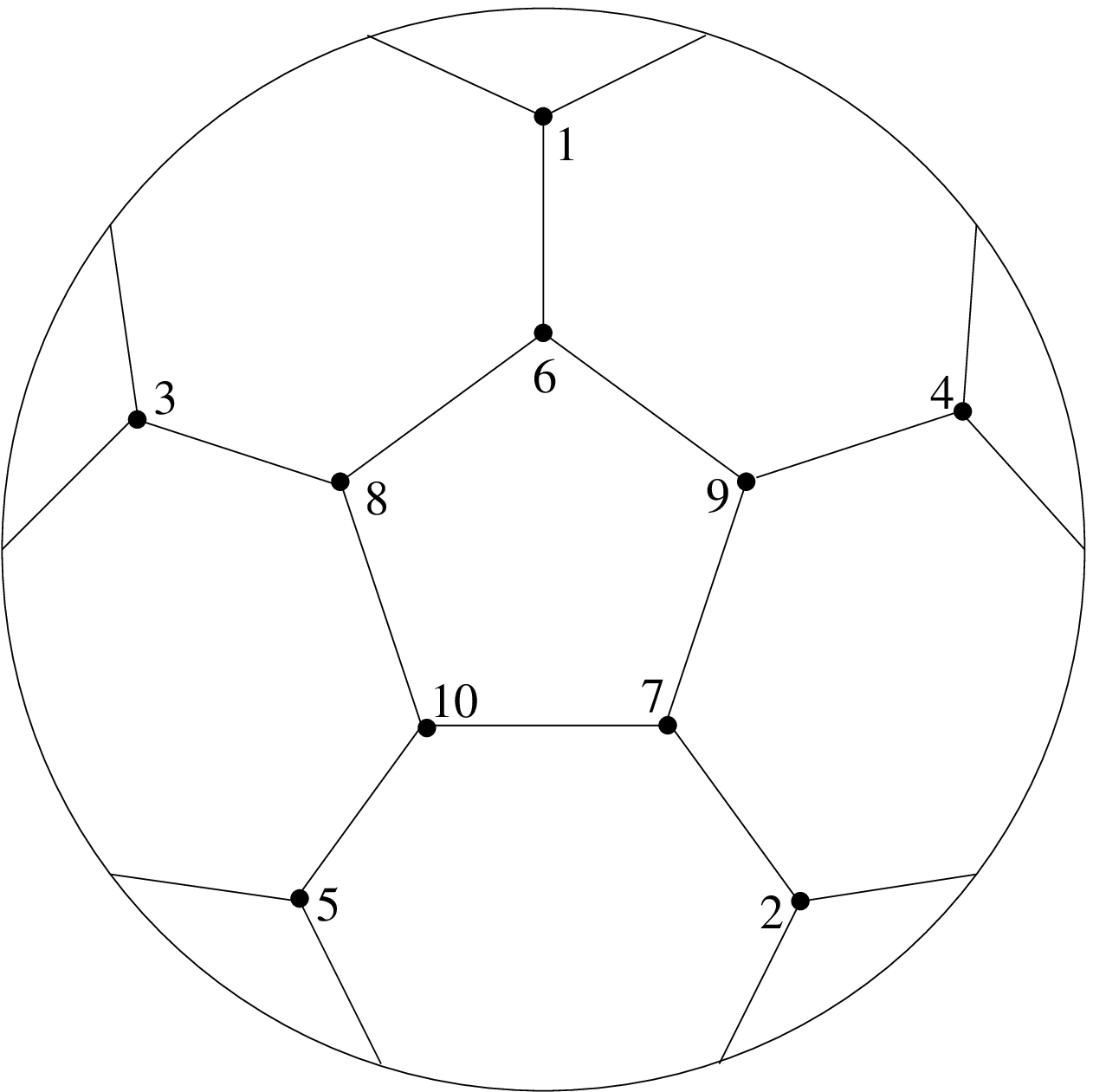}
  \caption{One of the two projective-planar embeddings of the Petersen graph}
  \label{fig:petersen}
\end{figure}

Let $P_{10}$ denote a labeling of the Petersen graph as shown in
\figref{fig:petersen}. In fact, \figref{fig:petersen} shows an
embedding of $P_{10}$ in the projective plane. The disk system \D\
associated with this embedding consists of the 5-cycles 6-9-7-10-8,
1-5-10-7-2, 4-3-8-10-5, 2-1-6-8-3,  5-4-9-6-1, and 3-2-7-9-4.

$P_{10}$ has exactly one other embedding in the projective plane. This
embedding is distinct from the above embedding, but is isomorphic to
it. (An isomorphism of embeddings is an isomorphism $\tau$ of the
underlying graphs such that a cycle $C$ is facial in one embedding if
and only if $\tau(C)$ is facial in the other.) The disk system $\D'$
associated with the second embedding consists of the 5-cycles
1-2-3-4-5, 6-9-4-3-8, 7-10-5-4-9, 8-6-1-5-10, 9-7-2-1-6, and
10-8-3-2-7.

We now describe notation that will let us denote specific enlargements
of a (labeled) graph as given by \thmref{thm:main}. Recall the
operations~\mbox{1--9} and the definition of a split, as described in
Sections~\ref{sec:disks-intro} and \ref{sec:main}.

Let $G$ be a graph whose vertices are labeled $1,\ldots,n$. For
vertices $u,v$, the graph $G\!+\!(u,v)$ denotes the graph obtained
from $G$ by adding an edge joining $u$ and $v$ (if none existed
before). Also, the graph $G*v(N_1)$ denotes the graph obtained by
splitting the vertex $v$, where $N_1$ is as in the definition of a
split. We follow the convention that the vertex $v_1$ retains the same
label as $v$, while $v_2$ is assigned the label $n+1$.

Since operations~1--7 are defined in terms of vertex splits and edge
additions, the above notation lets us specify \ext{$i$}s for
$i=1,\ldots,7$. An \ext{\triad} of $G$ is specified as
$G\!+\!(x_1,x_2,x_3)$, where the vertices $x_i$ are as in the
definition of operation~\triad. The new vertex $x$ gets the label
$n+1$.

Finally, a \ext{\tedge} of $G$ is specified as $G\!+\!(u,x\!-\!y)$,
where $u,x,y$ are as in the definition of operation~\tedge. The new
vertex obtained by subdividing the edge $xy$ gets the label $n+1$.

\subsection{Proof of Theorem~\ref{thm:petersen}.}

For the backward implication of \thmref{thm:petersen}, recall that
each of the eight graphs specified is either isomorphic to one of the
35 forbidden minors of \cite{arch-proj} or is obtained from one of
them by splitting a vertex. In particular, none of these eight graphs
embed in the projective plane, and so $H$ does not embed either.

For the forward implication, $H$, by hypothesis, does not embed in the
projective plane, and has a minor isomorphic to $P_{10}$. Clearly, the
disk system \D\ of $P_{10}$ has no locally planar extension to $H$.
Applying \thmref{thm:main} to $P_{10},\D$ and $H$, it is easy to check
that $H$ has a minor isomorphic to one of three enlargements, up to
isomorphism:
\begin{enumerate}
\item a \ext{\cross} $Q_1=P_{10}+(7,8)+(9,10)$
\item an \ext{\triad} $Q_2=P_{10}+(2,4,6)$
\item a \ext{\tedge} $Q_3=P_{10}+(1,3-4)$
\end{enumerate}

$Q_2$ has a minor isomorphic to \kffminus, as witnessed by the branch
sets $\{1,5\}$, $\{3,8\}$, $\{7,9\}$, $\{2\}$, $\{4\}$, $\{6\}$,
$\{10\}$, and $\{11\}$. (The order of the branch sets follows that of
the corresponding vertex labels in \kffminus, as shown in
\figref{fig:proj-list}.)

\begin{figure}[ht]
\centering
\includegraphics[width=0.4\textwidth]{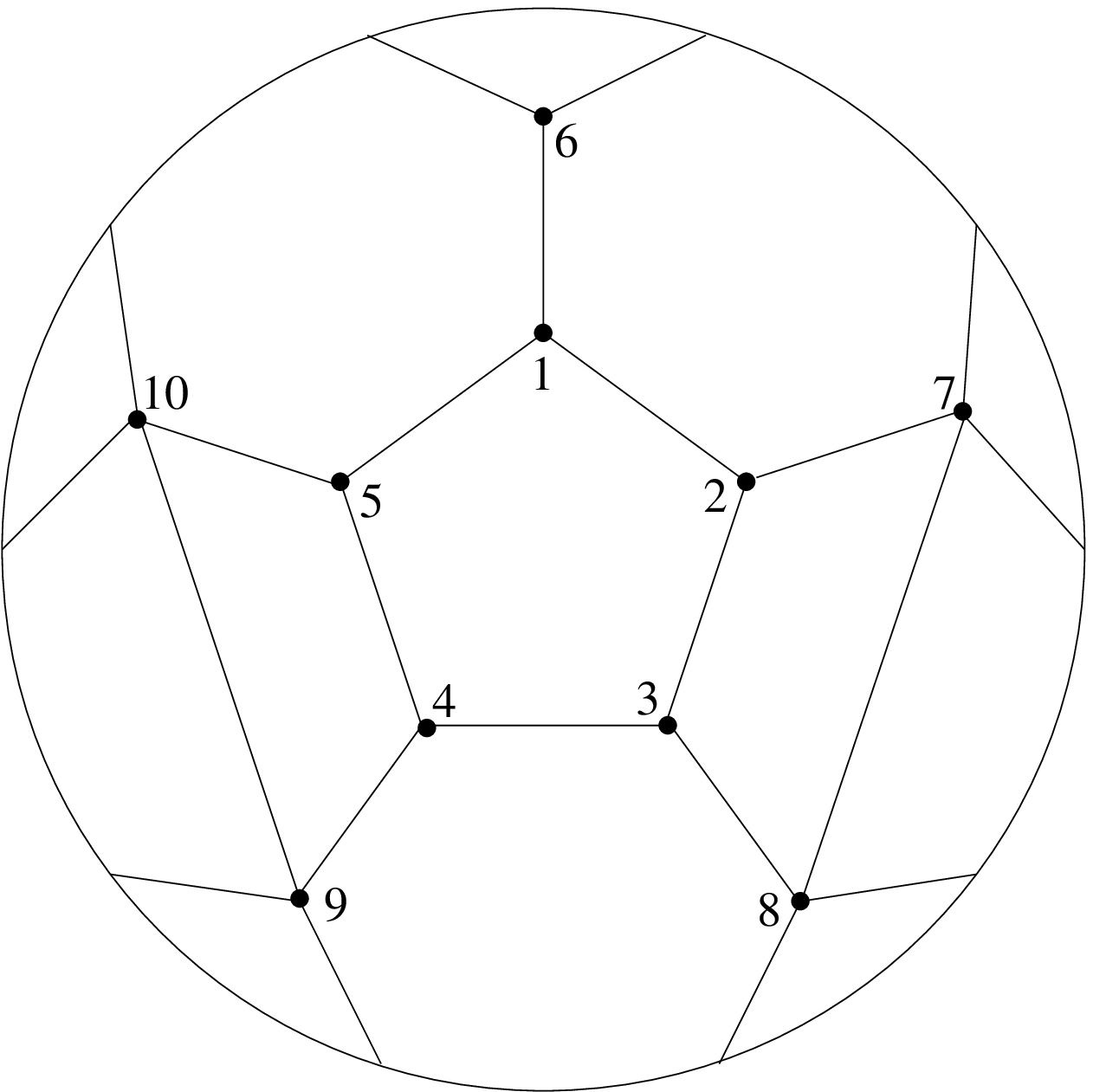}\hfill%
\includegraphics[width=0.4\textwidth]{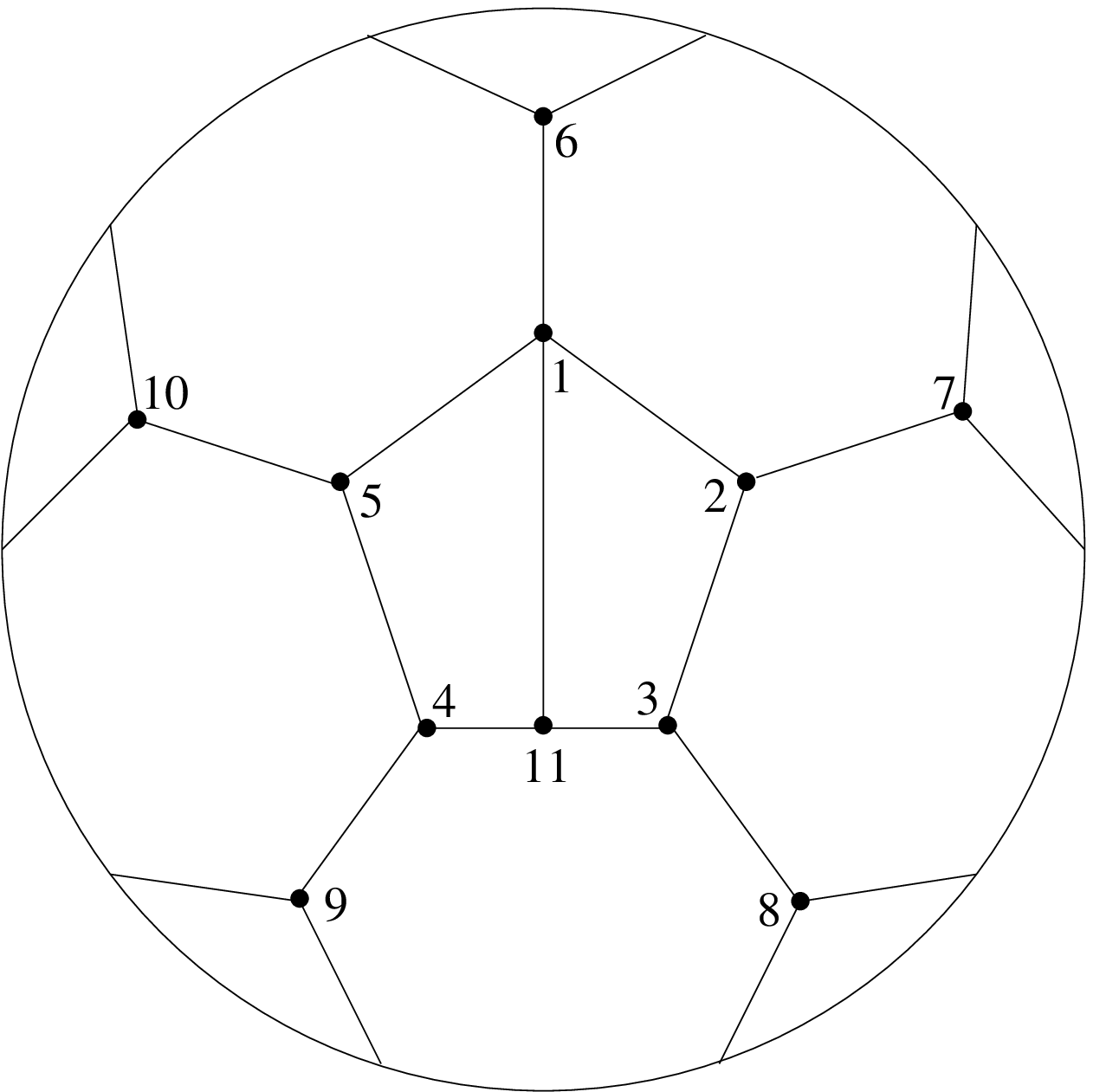}\\
\makebox[0.4\textwidth][c]{$Q_1$}\hfill%
\makebox[0.4\textwidth][c]{$Q_3$}
\caption{The graphs $Q_1$ and $Q_3$}
\label{fig:petersen-cross-tedge}
\end{figure}

Thus we may assume that $H$ has a minor isomorphic to $Q_1$ or
$Q_3$. The disk system $\D'$ of $P_{10}$ extends in a natural way to
disk systems $\D_1,\D_3$ in the enlargements $Q_1,Q_3$. Thus $Q_1,Q_3$
each embed (\emph{uniquely}) in the projective plane. The embeddings
are shown in \figref{fig:petersen-cross-tedge}.

We now apply \thmref{thm:main} to $Q_1,\D_1,H$ and $Q_3,\D_3,H$ and
deduce \thmref{thm:petersen}. This involves a fair amount of
case-checking, which is summarized in Tables~\ref{tab:q1} and
\ref{tab:q3}. Each row in the tables lists an enlargement of $Q_1$ or
$Q_3$, along with one of the eight graphs from the list that is a
minor of the enlargement. The branch sets in the rightmost column
follow the order of the vertex labels of the corresponding graph in
the preceding column. For clarity, singleton sets are not enclosed in
braces.

Tables~\ref{tab:q1} and \ref{tab:q3} respectively list all possible
enlargements of $Q_1$ and $Q_3$ up to isomorphism, with the exception
of \ext{\triad}s and \ext{\tedge}s of $Q_1$, and \ext{\triad}s of
$Q_3$. Every \ext{\triad} of $Q_1$ with respect to $\D_1$ has a
subgraph isomorphic to $Q_2$, and thus has a minor isomorphic to
\kffminus. Every \ext{\triad} of $Q_3$ with respect to $\D_3$ either
has a minor isomorphic to $Q_2$, or is isomorphic to the \ext{\triad}
listed in \tabref{tab:q3}. Finally, every \ext{\tedge} of $Q_1$ with
respect to $\D_1$ is either isomorphic to the \ext{\tedge} listed in
\tabref{tab:q1} or is isomorphic to a \ext{\cross} of $Q_3$ with
respect to $\D_3$ (and is thus listed in \tabref{tab:q3}
instead). This finishes the proof of \thmref{thm:petersen}.\qed
\newcommand{\myraise}[1]{\raisebox{1.2ex}[0in]{#1}}
\begin{table}[ht]
\centering
\caption{Applying \thmref{thm:main} to $Q_1$}
\label{tab:q1}
\begin{tabular}{|c|c|c|c|}
\hline
Type&Enlargement&Minor&Branch sets of the minor\\
\hline
& $Q_1\!+\!(2,10)$ & \dsplit & $\{2,3\},7,9,8,10,1,5,4,6$\\
\myraise{\jump}& $Q_1\!+\!(3,10)$ & \fsplit & $8,7,2,3,10,1,9,4,5,6$\\
\hline
& $Q_1\!+\!(2,8)\!+\!(3,7)$ & \etwenty & $2,7,3,8,\{1,6\},9,4,10,5$\\
& $Q_1\!+\!(2,4)\!+\!(3,5)$ & \kffmfour & $2,3,5,4,7,8,10,9,\{1,6\}$\\
& $Q_1\!+\!(1,4)\!+\!(3,5)$ & \ffour & $2,4,5,1,7,9,10,6,8,3$\\
& $Q_1\!+\!(1,4)\!+\!(2,5)$ & \ffour & $1,3,5,2,6,8,10,7,9,4$\\
& $Q_1\!+\!(3,9)\!+\!(4,8)$ & \cthree & $3,4,1,10,7,9,2,5,\{6,8\}$\\
& $Q_1\!+\!(3,9)\!+\!(4,6)$ & \cthree & $3,4,1,10,7,9,2,5,\{6,8\}$\\
\myraise{\cross}& $Q_1\!+\!(4,6)\!+\!(8,9)$ & \etwenty & $8,7,10,9,3,\{1,2\},5,6,4$\\
& $Q_1\!+\!(2,9)\!+\!(6,7)$ & \dsplit & $\{1,2\},7,10,6,9,3,4,5,8$\\
& $Q_1\!+\!(1,9)\!+\!(6,7)$ & \fsplit & $1,5,4,9,10,3,7,6,8,2$\\
& $Q_1\!+\!(1,9)\!+\!(2,6)$ & \ffour & $1,10,4,9,6,8,3,7,2,5$\\
& $Q_1\!+\!(1,7)\!+\!(2,9)$ & \dsplit & $9,7,8,2,\{1,6\},4,5,10,3$\\
& $Q_1\!+\!(1,7)\!+\!(2,6)$ & \cthree & $1,2,4,8,10,7,5,3,\{6,9\}$\\
\hline
& $Q_1\!*\!7(2,10)$ & \fsplit & $\{1,6\},5,4,9,10,3,7,11,8,2$\\
\myraise{\ncsplit}& $Q_1\!*\!8(3,10)$ & \fsplit & $2,7,11,\{1,6\},9,8,4,5,10,3$\\
\hline
& $Q_1\!*\!7(2,9)\!+\!(1,11)$ & \fsplit & $2,7,11,\{1,6\},9,8,4,5,10,3$\\
& $Q_1\!*\!7(2,9)\!+\!(6,11)$ & \fsplit & $3,4,5,\{8,10\},9,1,7,11,6,2$\\
& $Q_1\!*\!7(2,8)\!+\!(3,11)$ & \fsplit & $8,7,2,3,11,1,9,4,\{5,10\},6$\\
& $Q_1\!*\!8(3,7)\!+\!(2,11)$ & \fsplit & $\{5,10\},1,6,11,2,9,3,8,7,4$\\
\myraise{\splitjump}& $Q_1\!*\!8(3,7)\!+\!(1,8)$ & \fsplit & $\{5,10\},11,6,1,8,9,3,2,7,4$\\
& $Q_1\!*\!8(3,7)\!+\!(5,8)$ & \fsplit & $\{1,2\},3,4,5,8,9,11,10,7,6$\\
& $Q_1\!*\!8(3,6)\!+\!(4,11)$ & \fsplit & $8,3,\{2,7\},11,4,1,9,10,5,6$\\
& $Q_1\!*\!8(3,6)\!+\!(9,11)$ & \etwenty & $7,10,9,11,2,\{1,5,6\},4,8,3$\\
\hline
& $Q_1\!*\!7(8,10)\!*\!8(3,7)\!+\!(11,12)$ & \fsplit & $\{1,2\},6,9,11,12,\{3,4\},10,7,8,5$\\
& $Q_1\!*\!7(2,8)\!*\!8(7,10)\!+\!(11,12)$ & \fsplit & $\{3,4\},9,6,12,11,\{1,2\},10,8,7,5$\\
& $Q_1\!*\!7(2,9)\!*\!9(4,6)\!+\!(9,11)$ & \fsplit & $\{3,4\},8,6,9,11,\{1,2\},10,12,7,5$\\
\myraise{\dsplitjump}& $Q_1\!*\!7(2,8)\!*\!9(4,10)\!+\!(7,9)$ & \fsplit & $8,\{3,4\},5,10,9,\{1,2\},12,11,7,6$\\
& $Q_1\!*\!7(2,9)\!*\!10(9,11)\!+\!(7,12)$ & \fsplit & $\{1,6\},2,3,\{8,11\},7,4,12,10,9,5$\\
& $Q_1\!*\!7(2,8)\!*\!10(5,9)\!+\!(7,10)$ & \fsplit & $3,\{1,2\},6,8,7,9,10,12,11,\{4,5\}$\\
\hline
& $Q_1\!*\!7(2,8)\!+\!(1,11)\!+\!(6,7)$ & \fsplit & $2,7,6,1,11,8,\{4,9\},5,10,3$\\
\splitcross& $Q_1\!*\!8(3,7)\!+\!(4,11)\!+\!(8,9)$ & \fsplit & $\{1,6\},5,10,11,4,7,3,8,9,2$\\
& $Q_1\!*\!8(3,6)\!+\!(1,11)\!+\!(8,10)$ & \fsplit & $2,7,11,\{1,6\},9,8,4,5,10,3$\\
\hline
\dsplitcross& $Q_1\!*\!8(3,7)\!*\!9(4,10)\!+\!(8,12)\!+\!(9,11)$ & \fsplit & $\{2,7\},10,5,\{1,6\},11,4,8,12,9,3$\\
\hline
\tedge& $Q_1\!+\!(1,7\!-\!8)$ & \fsplit & $2,7,11,\{1,6\},9,8,4,5,10,3$\\
\hline
\end{tabular}
\end{table}

\begin{table}[ht]
\centering
\caption{Applying \thmref{thm:main} to $Q_3$}
\label{tab:q3}
\begin{tabular}{|c|c|c|c|}
\hline
Type&Enlargement&Minor&Branch sets of the minor\\
\hline
& $Q_3\!+\!(2,4)$ & \fsplit & $1,11,3,2,\{4,5\},8,9,7,10,6$\\
\myraise{\jump}& $Q_3\!+\!(2,5)$ & \fsplit & $1,11,4,5,\{2,3\},9,8,10,7,6$\\
\hline
& $Q_3\!+\!(1,7)\!+\!(2,6)$ & \dsplit & $\{3,8,11\},2,7,6,1,4,5,10,9$\\
& $Q_3\!+\!(1,9)\!+\!(2,6)$ & \fsplit & $9,4,5,1,\{3,11\},10,2,6,8,7$\\
& $Q_3\!+\!(1,9)\!+\!(6,7)$ & \ffour & $5,11,9,1,10,\{3,8\},6,2,7,4$\\
& $Q_3\!+\!(2,9)\!+\!(6,7)$ & \kffminus & $1,\{3,8\},\{4,9\},2,\{5,10\},6,11,7$\\
& $Q_3\!+\!(1,7)\!+\!(2,9)$ & \dsplit & $1,2,\{3,8\},7,\{6,9\},5,4,11,10$\\
& $Q_3\!+\!(2,8)\!+\!(3,7)$ & \ffour & $11,9,5,\{1,6\},3,7,10,2,8,4$\\
& $Q_3\!+\!(2,10)\!+\!(3,7)$ & \kffmfour & $1,2,3,11,5,10,7,\{4,9\},\{6,8\}$\\
& $Q_3\!+\!(2,10)\!+\!(7,8)$ & \ffour & $3,7,10,8,11,\{4,9\},5,6,1,2$\\
& $Q_3\!+\!(3,10)\!+\!(7,8)$ & \ffour & $5,11,6,1,10,3,8,2,7,\{4,9\}$\\
& $Q_3\!+\!(2,8)\!+\!(3,10)$ & \fsplit & $2,8,6,1,\{3,11\},9,10,5,4,7$\\
\myraise{\cross}& $Q_3\!+\!(3,9)\!+\!(4,8)$ & \kffmfour & $\{1,6\},2,3,11,5,\{7,10\},8,4,9$\\
& $Q_3\!+\!(3,9)\!+\!(8,11)$ & \fsplit & $11,4,5,\{1,6\},9,10,3,2,7,8$\\
& $Q_3\!+\!(3,4)\!+\!(8,11)$ & \dsplit & $\{4,5\},11,8,1,\{2,3\},9,7,10,6$\\
& $Q_3\!+\!(6,11)\!+\!(8,9)$ & \fsplit & $10,7,2,\{3,8\},9,1,4,11,6,5$\\
& $Q_3\!+\!(3,9)\!+\!(6,11)$ & \fsplit & $10,7,2,\{3,8\},9,1,4,11,6,5$\\
& $Q_3\!+\!(3,4)\!+\!(6,11)$ & \dsplit & $\{1,2\},11,4,6,\{3,8\},7,10,5,9$\\
& $Q_3\!+\!(4,6)\!+\!(8,9)$ & \ffour & $5,11,6,1,10,\{3,8\},9,2,7,4$\\
& $Q_3\!+\!(3,9)\!+\!(4,6)$ & \ffour & $5,11,6,1,10,\{3,8\},9,2,7,4$\\
& $Q_3\!+\!(3,6)\!+\!(8,11)$ & \etwenty & $2,1,3,11,7,\{6,9\},8,\{4,5\},10$\\
& $Q_3\!+\!(1,3)\!+\!(2,11)$ & \fsplit & $2,3,6,1,11,9,\{8,10\},5,4,7$\\
\hline
\ncsplit& $Q_3\!*\!1(2,5)$ & \fsplit & $7,10,5,\{1,2\},\{3,8\},4,6,12,11,9$\\
\hline
& $Q_3\!*\!1(5,6)\!+\!(10,12)$ & \fsplit & $12,1,5,10,\{6,8\},4,\{2,3\},7,9,11$\\
\splitjump& $Q_3\!*\!1(5,6)\!+\!(8,12)$ & \fsplit & $9,\{1,6\},5,\{4,11\},12,10,2,3,8,7$\\
& $Q_3\!*\!1(5,6)\!+\!(1,3)$ & \fsplit & $10,8,6,\{1,5\},3,\{4,9\},2,12,11,7$\\
\hline
\splitcross& $Q_3\!*\!1(5,6)\!+\!(1,7)\!+\!(9,12)$ & \fsplit & $8,10,5,\{1,6\},7,4,2,12,9,\{3,11\}$\\
\hline
\triad& $Q_3+(2,9,11)$ & \ffour & $1,12,3,2,\{4,5\},9,\{6,8\},7,10,11$\\
\hline
& $Q_3\!+\!(8,1\!-\!11)$ & \fsplit & $1,12,11,\{2,3\},\{6,8\},4,10,7,9,5$\\
& $Q_3\!+\!(8,1\!-\!2)$ & \fsplit & $9,4,5,\{1,6\},\{3,11\},10,2,12,8,7$\\
& $Q_3\!+\!(10,1\!-\!2)$ & \ffour & $11,5,9,\{1,6\},\{3,8\},10,7,12,2,4$\\
& $Q_3\!+\!(6,2\!-\!3)$ & \fsplit & $2,12,6,1,\{3,8,11\},9,10,5,4,7$\\
& $Q_3\!+\!(9,2\!-\!3)$ & \fsplit & $5,4,11,\{1,6\},9,\{3,8\},7,2,12,10$\\
& $Q_3\!+\!(3,1\!-\!6)$ & \ffour & $2,11,12,1,7,\{4,9\},\{6,8\},5,10,3$\\
\tedge& $Q_3\!+\!(1,3\!-\!8)$ & \ffour & $2,11,12,1,7,\{4,6,9\},8,5,10,3$\\
& $Q_3\!+\!(2,6\!-\!8)$ & \ffour & $3,7,12,\{8,10\},11,\{4,9\},6,5,1,2$\\
& $Q_3\!+\!(7,6\!-\!8)$ & \ffour & $11,9,5,\{1,6\},\{2,3\},7,10,12,8,4$\\
& $Q_3\!+\!(3,7\!-\!9)$ & \ffour & $5,11,9,\{1,6\},10,\{3,8\},12,2,7,4$\\
& $Q_3\!+\!(8,7\!-\!9)$ & \ffour & $5,11,9,\{1,6\},10,\{3,8\},12,2,7,4$\\
& $Q_3\!+\!(1,7\!-\!10)$ & \etwo & $2,9,12,11,5,8,\{1,6\},3,7,4,10$\\
& $Q_3\!+\!(6,7\!-\!10)$ & \etwo & $2,9,12,11,5,8,\{1,6\},3,7,4,10$\\
\hline
\end{tabular}
\end{table}
\clearpage

\section*{Acknowledgment}
The results of this paper form part of the Ph.D.\ dissertation~\cite{HegdePhD}
written by the first author under the direction of the second author.

\bibliographystyle{plain}
\nocite{*}

\begin{thebibliography}{10}

\bibitem{arch-phd}
D. Archdeacon.
\newblock {\em A {K}uratowski Theorem for the Projective Plane}.
\newblock PhD thesis, The {O}hio {S}tate {U}niversity, 1980.

\bibitem{arch-proj}
D. Archdeacon.
\newblock A {K}uratowski theorem for the projective plane.
\newblock {\em J. Graph Theory}, 5(3):243--246, 1981.

\bibitem{DinOpoThoVer}
G.~Ding, B.~Oporowski, R.~Thomas, and D.~Vertigan.
\newblock Large non-planar graphs and an application to crossing-critical graphs.
\newblock {\it  J.~Combin.\ Theory Ser.\ B} {\bf 101} (2011), 111--121.

\bibitem{HegdePhD}
R. Hegde.
\newblock {\em New Tools and Results in Graph Structure Theory}.
\newblock PhD thesis, Georgia {I}nstitute of {T}echnology, 2006.

\bibitem{JohThoSplit}
T.~Johnson and R.~Thomas.
\newblock A splitter theorem for internally $4$-connected graphs.
\newblock Manuscript.

\bibitem{KawNorThoWolbdtw}
K.~Kawarabayashi, S.~Norine, R.~Thomas and P.~Wollan.
\newblock $K_6$ minors in $6$-connected graphs of bounded tree-width.
\newblock {\tt arXiv:1203.2171}.

\bibitem{RobSeyThoExt}
S.~Norin and R.~Thomas.
\newblock Non-planar extensions of planar graphs.
\newblock Manuscript.

\bibitem{ThoThoTutte}
R.~Thomas and J.~Thomson.
\newblock Excluding minors in non-planar graphs of girth five.
\newblock {\em Combinatorics, Probability and Computing}, 9:573--585, 2000.

\bibitem{TutHowto}
W.~T. Tutte.
\newblock How to draw a graph.
\newblock {\em Proc. London Math. Soc. (3)}, 13:743--767, 1963.

\bibitem{WhiCongr}
H.~Whitney.
\newblock Congruent graphs and the connectivity of graphs.
\newblock {\em Amer. J. Math.}, 54(1):150--168, 1932.

\end{thebibliography}

\baselineskip 11pt
\vfill
\noindent
This material is based upon work supported by the National Science Foundation
under Grants No.~0200595 and~0354742. 
Any opinions, findings, and conclusions or
recommendations expressed in this material are those of the authors and do
not necessarily reflect the views of the National Science Foundation.
\eject

\end{document}